\documentclass[10pt]{article}
\usepackage[centertags]{amsmath}
\usepackage{array,amsmath}
\usepackage{amsfonts}
\usepackage{amssymb}
\usepackage{amsthm}
\usepackage{newlfont}
\usepackage{hyperref}
\usepackage[dvips]{graphicx}
\usepackage{pst-node}
\usepackage{tikz}
\usepackage[all]{xy}
\usepackage{color}
\usepackage{mathdots}
\newtheorem{thm}{Theorem}[section]
\newtheorem{pro}[thm]{Proposition}
\newtheorem{lem}[thm]{Lemma}
\newtheorem{rem}[thm]{Remark}
\newtheorem{rems}[thm]{Remarks}
\newtheorem{cor}[thm]{Corollary}
\newtheorem{df}[thm]{Definition}
\newtheorem*{prf}{Proof}
\newtheorem{ex}[thm]{Example}
\newtheorem{exs}[thm]{Examples}

 \numberwithin{equation}{section}

 \renewcommand{\a }{\alpha }
\renewcommand{\b }{\beta }

\newcommand{\ld }{\lambda }

\renewcommand{\t }{\theta}
\newcommand{\Z}{\mathbb{Z}}
\newcommand{\K}{\mathbb{K}}

\newcommand{\N}{\mathbb{N}}

\newcommand{\tone}{{\begin{tikzpicture}[ scale=0.05]
\fill[color=black](1, 0)circle(0mm);
\fill[color=black](1, 10)circle(0mm);
\draw[-](1, 0)--(1, 10);
\end{tikzpicture}}}

\newcommand{\ttwo}{\begin{tikzpicture}[ scale=0.05]
\fill[color=black](8, 0)circle(0mm);
\fill[color=black](8, 6)circle(0mm);
\fill[color=black](0, 14)circle(0mm);
\fill[color=black](16, 14)circle(0mm);
\draw[-](8, 0)--(8, 6);
\draw[-](8, 6)--(0, 14);
\draw[-](8, 6)--(16, 14);
\end{tikzpicture}}

\newcommand{\tthreeone}{\begin{tikzpicture}[ scale=0.05]
\fill[color=black](8, 0)circle(0mm);
\fill[color=black](8, 6)circle(0mm);
\fill[color=black](0, 14)circle(0mm);
\fill[color=black](16, 14)circle(0mm);
\fill[color=black](12, 10)circle(0mm);
\fill[color=black](8, 14)circle(0mm);
\draw[-](8, 0)--(8, 6);
\draw[-](8, 6)--(0, 14);
\draw[-](8, 6)--(16, 14);
\draw[-](12, 10)--(8, 14);
\end{tikzpicture}}

\newcommand{\tthreetwo}{\begin{tikzpicture}[ scale=0.05]
\fill[color=black](8, 0)circle(0mm);
\fill[color=black](8, 6)circle(0mm);
\fill[color=black](0, 14)circle(0mm);
\fill[color=black](16, 14)circle(0mm);
\fill[color=black](4, 10)circle(0mm);
\fill[color=black](8, 14)circle(0mm);
\draw[-](8, 0)--(8, 6);
\draw[-](8, 6)--(0, 14);
\draw[-](4, 10)--(8, 14);
\draw[-](8, 6)--(16, 14);
\end{tikzpicture}}

\newcommand{\tfourone}{\begin{tikzpicture}[ scale=0.05]
\fill[color=black](9, -4)circle(0mm);
\fill[color=black](9, 2)circle(0mm);
\fill[color=black](0, 11)circle(0mm);
\fill[color=black](18, 11)circle(0mm);
\fill[color=black](12, 5)circle(0mm);
\fill[color=black](6, 11)circle(0mm);
\fill[color=black](15, 8)circle(0mm);
\fill[color=black](12, 11)circle(0mm);
\draw[-](9, -4)--(9, 2);
\draw[-](9, 2)--(0, 11);
\draw[-](9, 2)--(18, 11);
\draw[-](12, 5)--(6, 11);
\draw[-](15, 8)--(12, 11);
\end{tikzpicture}}

\newcommand{\tfourtwo}{\begin{tikzpicture}[ scale=0.05]
\fill[color=black](9, -4)circle(0mm);
\fill[color=black](9, 2)circle(0mm);
\fill[color=black](0, 11)circle(0mm);
\fill[color=black](18, 11)circle(0mm);
\fill[color=black](12, 5)circle(0mm);
\fill[color=black](6, 11)circle(0mm);
\fill[color=black](9, 8)circle(0mm);
\fill[color=black](12, 11)circle(0mm);
\draw[-](9, -4)--(9, 2);
\draw[-](9, 2)--(0, 11);
\draw[-](9, 2)--(18, 11);
\draw[-](12, 5)--(6, 11);
\draw[-](9, 8)--(12, 11);
\end{tikzpicture}}

\newcommand{\tfourthree}{\begin{tikzpicture}[ scale=0.05]
\fill[color=black](9, -4)circle(0mm);
\fill[color=black](9, 2)circle(0mm);
\fill[color=black](0, 11)circle(0mm);
\fill[color=black](3, 8)circle(0mm);
\fill[color=black](6, 11)circle(0mm);
\fill[color=black](18, 11)circle(0mm);
\fill[color=black](15, 8)circle(0mm);
\fill[color=black](12, 11)circle(0mm);
\draw[-](9, -4)--(9, 2);
\draw[-](9, 2)--(0, 11);
\draw[-](3, 8)--(6, 11);
\draw[-](9, 2)--(18, 11);
\draw[-](15, 8)--(12, 11);
\end{tikzpicture}}

\newcommand{\tfourfour}{\begin{tikzpicture}[ scale=0.05]
\fill[color=black](9, -4)circle(0mm);
\fill[color=black](9, 2)circle(0mm);
\fill[color=black](0, 11)circle(0mm);
\fill[color=black](6, 5)circle(0mm);
\fill[color=black](12, 11)circle(0mm);
\fill[color=black](9, 8)circle(0mm);
\fill[color=black](6, 11)circle(0mm);
\fill[color=black](18, 11)circle(0mm);
\draw[-](9, -4)--(9, 2);
\draw[-](9, 2)--(0, 11);
\draw[-](6, 5)--(12, 11);
\draw[-](9, 8)--(6, 11);
\draw[-](9, 2)--(18, 11);
\end{tikzpicture}}

\newcommand{\tfourfive}{\begin{tikzpicture}[ scale=0.05]
\fill[color=black](9, -4)circle(0mm);
\fill[color=black](9, 2)circle(0mm);
\fill[color=black](0, 11)circle(0mm);
\fill[color=black](3, 8)circle(0mm);
\fill[color=black](6, 11)circle(0mm);
\fill[color=black](6, 5)circle(0mm);
\fill[color=black](12, 11)circle(0mm);
\fill[color=black](18, 11)circle(0mm);
\draw[-](9, -4)--(9, 2);
\draw[-](9, 2)--(0, 11);
\draw[-](3, 8)--(6, 11);
\draw[-](6, 5)--(12, 11);
\draw[-](9, 2)--(18, 11);
\end{tikzpicture}}

\newcommand{\psiwedgephi}{\begin{tikzpicture}[ scale=0.05]
\fill[color=black](8, 0)circle(0mm);
\fill[color=black](8, 6)circle(0mm);
\fill[color=black](0, 14)circle(0mm);
\fill[color=black](16, 14)circle(0mm);
\draw[-](8, 0)--(8, 6);
\draw[-](8, 6)--(0, 14);
\draw[-](8, 6)--(16, 14);
\path(-1,18) node {$\psi$};
\path(19,18) node {$\varphi$};
\end{tikzpicture}}
\newcommand{\ntree}{\begin{tikzpicture}[ scale=0.05]
\fill[color=black](8, 0)circle(0mm);
\fill[color=black](8, 6)circle(0mm);
\fill[color=black](0, 14)circle(0mm);
\fill[color=black](16, 14)circle(0mm);
\draw[-](8, 0)--(8, 6);
\draw[-](8, 6)--(0, 14);
\draw[-](8, 6)--(16, 14);
\path(-8,18) node {\small{L-subtree}};
\path(25,18) node {\small{R-subtree}};
\end{tikzpicture}}
\newcommand{\LAL}{\begin{tikzpicture}[ scale=0.1]
\fill[color=black](9, -4)circle(0mm);
\fill[color=black](9, 2)circle(0mm);
\fill[color=black](0, 11)circle(0mm);
\fill[color=black](6, 5)circle(0mm);
\fill[color=black](12, 11)circle(0mm);
\fill[color=red](9, 8)circle(5mm);
\fill[color=black](6, 11)circle(0mm);
\fill[color=black](18, 11)circle(0mm);
\draw[dashed](9, -4)--(9, 2);
\draw[dashed](9, 2)--(6, 5);
\draw[dashed](6, 5)--(0, 11);
\draw[dashed](6, 5)--(9, 8);
\draw[-, color=blue](9, 8)--(12, 11);
\draw[-, color=blue](9, 8)--(6, 11);
\draw[dashed](9, 2)--(18, 11);
\path(9,14) node {\textcolor{blue}{$\overbrace{}^{LAL}$}};
\end{tikzpicture}
}
\newcommand{\RAL}{\begin{tikzpicture}[ scale=0.1]
\fill[color=black](9, -4)circle(0mm);
\fill[color=black](9, 2)circle(0mm);
\fill[color=black](0, 11)circle(0mm);
\fill[color=black](18, 11)circle(0mm);
\fill[color=black](12, 5)circle(0mm);
\fill[color=black](6, 11)circle(0mm);
\fill[color=red](9, 8)circle(5mm);
\fill[color=black](12, 11)circle(0mm);
\draw[dashed](9, -4)--(9, 2);
\draw[dashed](9, 2)--(0, 11);
\draw[dashed](9, 2)--(12, 5);
\draw[dashed](12, 5)--(18, 11);
\draw[dashed](12, 5)--(9, 8);
\draw[-,color=blue](9, 8)--(6, 11);
\draw[-, color=blue](9, 8)--(12, 11);
\path(9,14) node {\textcolor{blue}{$\overbrace{}^{RAL}$}};
\end{tikzpicture}}
\newcommand{\LDL}{\begin{tikzpicture}[ scale=0.1]
\fill[color=black](9, -4)circle(0mm);
\fill[color=black](9, 2)circle(0mm);
\fill[color=black](0, 11)circle(0mm);
\fill[color=red](3, 8)circle(5mm);
\fill[color=black](6, 11)circle(0mm);
\fill[color=red](6, 5)circle(5mm);
\fill[color=black](12, 11)circle(0mm);
\fill[color=black](18, 11)circle(0mm);
\draw[dashed](9, -4)--(9, 2);
\draw[dashed](9, 2)--(6, 5);
\draw[dashed](6, 5)--(3, 8);
\draw[dashed](3, 8)--(0, 11);
\draw[-,color=blue](3, 8)--(6, 11);
\draw[-,color=blue](6, 5)--(12, 11);
\draw[dashed](9, 2)--(18, 11);
\path(9,14) node {\textcolor{blue}{$\overbrace{}^{LDL}$}};
\end{tikzpicture}}
\newcommand{\RDL}{\begin{tikzpicture}[ scale=0.1]
\fill[color=black](9, -4)circle(0mm);
\fill[color=black](9, 2)circle(0mm);
\fill[color=black](0, 11)circle(0mm);
\fill[color=black](18, 11)circle(0mm);
\fill[color=red](12, 5)circle(5mm);
\fill[color=black](6, 11)circle(0mm);
\fill[color=red](15, 8)circle(5mm);
\fill[color=black](12, 11)circle(0mm);
\draw[dashed](9, -4)--(9, 2);
\draw[dashed](9, 2)--(0, 11);
\draw[dashed](9, 2)--(12, 5);
\draw[dashed](12, 5)--(15, 8);
\draw[dashed](15, 8)--(18, 11);
\draw[-,color=blue](12, 5)--(6, 11);
\draw[-,color=blue](15, 8)--(12, 11);
\path(9,14) node {\textcolor{blue}{$\overbrace{}^{RDL}$}};
\end{tikzpicture}}
\newcommand{\DLon}{\begin{tikzpicture}[ scale=0.1]
\fill[color=black](9, -4)circle(0mm);
\fill[color=black](9, 2)circle(0mm);
\fill[color=black](3, 8)circle(0mm);
\fill[color=black](6, 11)circle(0mm);
\fill[color=black](15, 8)circle(0mm);
\fill[color=black](12, 11)circle(0mm);
\fill[color=black](0, 11)circle(0mm);
\fill[color=black](18, 11)circle(0mm);
\draw[dashed](9, -4)--(9, 2);
\draw[dashed](9, 2)--(3,8);
\draw[dashed](9, 2)--(15, 8);
\draw[-,color=blue](15,8)--(12,11);
\draw[dashed](15,8)--(18,11);
\draw[-,color=blue](3, 8)--(6, 11);
\draw[dashed](3, 8)--(0,11);
\path(9,14) node {\textcolor{blue}{$\overbrace{}^{DL1}$}};
\end{tikzpicture}}
\newcommand{\DLtw}{
\begin{tikzpicture}[ scale=0.1]
\fill[color=black](9, -4)circle(0mm);
\fill[color=black](9, 2)circle(0mm);
\fill[color=black](0, 11)circle(0mm);
\fill[color=black](4, 11)circle(0mm);
\fill[color=black](8, 11)circle(0mm);
\fill[color=black](12, 11)circle(0mm);
\fill[color=black](18, 11)circle(0mm);
\fill[color=black](4.5, 6)circle(0mm);
\fill[color=black](6, 8)circle(0mm);
\fill[color=black](15, 8)circle(1mm);
\draw[-](9, -4)--(9, 2);
\draw[dashed](9, 2)--(0, 11);
\draw[dashed](9, 2)--(18, 11);
\draw[-](4.5, 6)--(6, 8);
\draw[-](6, 8)--(4, 11) ;
\draw[-,color=blue](6, 8)--(8, 11);
\draw[-,color=blue](15, 8)--(12, 11);
\path(10,14) node {\textcolor{blue}{$\overbrace{}^{DL2}$}};
\path(18,12) node {$\iddots$};
\path(0,12) node {$\ddots$};
\end{tikzpicture}}
\newcommand{\DLtr}{
\begin{tikzpicture}[ scale=0.1]
\fill[color=black](9, -4)circle(0mm);
\fill[color=black](9, 2)circle(0mm);
\fill[color=black](0, 11)circle(0mm);
\fill[color=black](6, 11)circle(0mm);
\fill[color=black](10.6, 11)circle(0mm);
\fill[color=black](14, 11)circle(0mm);
\fill[color=black](18, 11)circle(0mm);
\fill[color=black](3, 8)circle(0mm);
\fill[color=black](13, 6)circle(0mm);
\fill[color=black](12, 8)circle(1mm);
\draw[-](9, -4)--(9, 2);
\draw[-](9, 2)--(3, 8);
\draw[-,color=blue](3, 8)--(6, 11);
\draw[dashed](13, 6)--(18, 11);
\draw[-](9, 2)--(13, 6);
\draw[-](13, 6)--(12, 8);
\draw[-,color=blue](12, 8)--(10.6, 11);
\draw[-](12, 8)--(14, 11);
\path(8,14) node {\textcolor{blue}{$\overbrace{}^{DL3}$}};
\path(18,12) node {$\iddots$};
\path(0,12) node {$\ddots$};
\end{tikzpicture}
}
\newcommand{\DLfr}{
\begin{tikzpicture}[ scale=0.1]
\fill[color=black](9, -4)circle(0mm);
\fill[color=black](9, 2)circle(0mm);
\fill[color=black](0, 11)circle(0mm);
\fill[color=black](6, 8)circle(0mm);
\fill[color=black](8, 11)circle(0mm);
\fill[color=black](4, 11)circle(0mm);
\fill[color=black](10.6, 11)circle(0mm);
\fill[color=black](14, 11)circle(0mm);
\fill[color=black](18, 11)circle(0mm);
\fill[color=black](5,6)circle(0mm);
\fill[color=black](13, 6)circle(0mm);
\fill[color=black](12, 8)circle(0mm);
\draw[-](9, -4)--(9, 2);
\draw[-](9, 2)--(3,8);
\draw[dashed](3,8)--(0,11);
\draw[-](5,6)--(6,8);
\draw[-](6,8)--(4, 11);
\draw[-,color=blue](6,8)--(8, 11);
\draw[-](9, 2)--(13, 6);
\draw[dashed](13,6)--(18, 11);
\draw[-](13, 6)--(12, 8);
\draw[-,color=blue](12, 8)--(10.6, 11);
\draw[-](12, 8)--(14, 11);
\path(9,14) node {\textcolor{blue}{$\overbrace{}^{DL4}$}};
\end{tikzpicture}
}
\newcommand{\exttwoD}{
\begin{tikzpicture}[ scale=0.1]
\tikzstyle{terminal}=[circle,draw,thick,scale=0.15mm]
\tikzstyle{cercle}=[circle,draw,thick,fill=black, scale=0.15mm]
\fill[color=black](8, 0)circle(0mm);
\fill[color=black](8, 6)circle(0mm);
\node[cercle] (C) at (16,15){};
\node[terminal] (T) at (1,15){};
\draw[-](8, 0)--(8, 6);
\draw[-](8, 6)--(1, 14);
\draw[-](8, 6)--(16, 14);
\end{tikzpicture}}
\newcommand{\extthreeD}{
\begin{tikzpicture}[ scale=0.1]
\tikzstyle{terminal}=[circle,draw,thick,scale=0.15mm]
\tikzstyle{cercle}=[circle,draw,thick,fill=black, scale=0.15mm]
\fill[color=black](8, 0)circle(0mm);
\fill[color=black](8, 6)circle(0mm);
\fill[color=black](4, 10)circle(0mm);
\node[cercle] (C) at (0,15){};
\node[cercle] (C) at (8,15){};
\node[terminal] (T) at (16, 15){};
\draw[-](8, 0)--(8, 6);
\draw[-](8, 6)--(0, 14);
\draw[-](4, 10)--(8, 14);
\draw[-](8, 6)--(16, 14);
\end{tikzpicture}
}
\newcommand{\extfiveD}{
\begin{tikzpicture}[ scale=0.1]
\tikzstyle{terminal}=[circle,draw,thick,scale=0.15mm]
\tikzstyle{cercle}=[circle,draw,thick,fill=black, scale=0.15mm]
\fill[color=black](9, -4)circle(0mm);
\fill[color=black](9, 2)circle(0mm);
\fill[color=black](4.5, 6)circle(0mm);
\fill[color=black](6, 8)circle(0mm);
\fill[color=black](15, 8)circle(0mm);
\node[cercle] (C) at (8,12){};
\node[terminal] (T) at (12,12){};
\node[cercle] (C) at (18,12){};
\node[terminal] (T) at (0, 12){};
\node[terminal] (T) at (4, 12){};
\draw[-](9, -4)--(9, 2);
\draw[-](9, 2)--(0, 11) ;
\draw[-](9, 2)--(18, 11);
\draw[-](4.5, 6)--(6, 8);
\draw[-](6, 8)--(4, 11);
\draw[-](6, 8)--(8, 11) ;
\draw[-](15, 8)--(12, 11);
\end{tikzpicture}}
\newcommand{\redoneone}{
\begin{tikzpicture}[ scale=0.05]
\tikzstyle{terminal}=[circle,draw,thick,scale=0.15mm]
\tikzstyle{cercle}=[circle,draw,thick,fill=black, scale=0.15mm]
\fill[color=black](1, 0)circle(0mm);
\node[cercle] (C) at (1,11){};
\draw[-](1, 0)--(1, 10);
\path (1,11)node[above]{$a$};
\end{tikzpicture}}
\newcommand{\redonetw}{
\begin{tikzpicture}[ scale=0.05]
\tikzstyle{terminal}=[circle,draw,thick,scale=0.15mm]
\tikzstyle{cercle}=[circle,draw,thick,fill=black, scale=0.15mm]
\fill[color=black](1, 0)circle(0mm);
\node[terminal] (T) at (1,11){};
\draw[-](1, 0)--(1, 10);
\path (1,11)node[above]{$a$};
\end{tikzpicture}}
\newcommand{\redonetwone}{
\begin{tikzpicture}[ scale=0.05]
\tikzstyle{terminal}=[circle,draw,thick,scale=0.15mm]
\tikzstyle{cercle}=[circle,draw,thick,fill=black, scale=0.15mm]
\fill[color=black](1, 0)circle(0mm);
\node[cercle] (C) at (1,11){};
\draw[-](1, 0)--(1, 10);
\path (1,11)node[above]{$a+1$};
\end{tikzpicture}}
\newcommand{\redtwtw}{
\begin{tikzpicture}[ scale=0.05]
\tikzstyle{terminal}=[circle,draw,thick,scale=0.15mm]
\tikzstyle{cercle}=[circle,draw,thick,fill=black, scale=0.15mm]
\fill[color=black](1, 0)circle(0mm);
\node[terminal] (T) at (1,11){};
\draw[-](1, 0)--(1, 10);
\path (1,11)node[above]{$a+1$};
\end{tikzpicture}}
\newcommand{\exttwo}{
\begin{tikzpicture}[ scale=0.05]
\tikzstyle{terminal}=[circle,draw,thick,scale=0.15mm]
\tikzstyle{cercle}=[circle,draw,thick,fill=black, scale=0.15mm]
\fill[color=black](8, 0)circle(0mm);
\fill[color=black](8, 6)circle(0mm);
\node[cercle] (C) at (16,15){};
\node[terminal] (T) at (0, 15){};
\draw[-](8, 0)--(8, 6);
\draw[-](8, 6)--(0, 14);
\draw[-](8, 6)--(16, 14);
\path (0,15) node[above]{$2$};
\path (16,15) node[above]{$3$};
\end{tikzpicture}}
\newcommand{\exttwotwo}{
\begin{tikzpicture}[ scale=0.05]
\tikzstyle{terminal}=[circle,draw,thick,scale=0.15mm]
\tikzstyle{cercle}=[circle,draw,thick,fill=black, scale=0.15mm]
\fill[color=black](8, 0)circle(0mm);
\fill[color=black](8, 6)circle(0mm);
\node[cercle] (C) at (16,15){};
\node[terminal] (T) at (0, 15){};
\draw[-](8, 0)--(8, 6);
\draw[-](8, 6)--(0, 14);
\draw[-](8, 6)--(16, 14);
\path (0,15) node[above]{$1$};
\path (16,15) node[above]{$2$};
\end{tikzpicture}}
\newcommand{\extredtwon}{
\begin{tikzpicture}[ scale=0.03]
\tikzstyle{terminal}=[circle,draw,thick,scale=0.15mm]
\tikzstyle{cercle}=[circle,draw,thick,fill=black, scale=0.15mm]
\fill[color=black](8, 0)circle(0mm);
\fill[color=black](8, 6)circle(0mm);
\node[cercle] (C) at (16,15){};
\node[terminal] (T) at (0, 15){};
\draw[-](8, 0)--(8, 6);
\draw[-](8, 6)--(0, 14);
\draw[-](8, 6)--(16, 14);
\path (0,15) node[above]{$a$};
\path (16,15) node[above]{$a$};
\end{tikzpicture}}
\newcommand{\extredttwon}{
\begin{tikzpicture}[ scale=0.03]
\tikzstyle{terminal}=[circle,draw,thick,scale=0.15mm]
\tikzstyle{cercle}=[circle,draw,thick,fill=black, scale=0.15mm]
\fill[color=black](8, 0)circle(0mm);
\fill[color=black](8, 6)circle(0mm);
\node[terminal] (T) at (16,15){};
\node[cercle] (C) at (0, 15){};
\draw[-](8, 0)--(8, 6);
\draw[-](8, 6)--(0, 14);
\draw[-](8, 6)--(16, 14);
\path (0,15) node[above]{$a$};
\path (16,15) node[above]{$a$};
\end{tikzpicture}}
\newcommand{\extredtwtw}{
\begin{tikzpicture}[ scale=0.03]
\tikzstyle{terminal}=[circle,draw,thick,scale=0.15mm]
\tikzstyle{cercle}=[circle,draw,thick,fill=black, scale=0.15mm]
\fill[color=black](8, 0)circle(0mm);
\fill[color=black](8, 6)circle(0mm);
\node[cercle] (C) at (16,15){};
\node[terminal] (T) at (0, 15){};
\draw[-](8, 0)--(8, 6);
\draw[-](8, 6)--(0, 14);
\draw[-](8, 6)--(16, 14);
\path (0,15) node[above]{$b$};
\path (16,15) node[above]{$a$};
\end{tikzpicture}}
\newcommand{\extredtwtwthr}{
\begin{tikzpicture}[ scale=0.03]
\tikzstyle{terminal}=[circle,draw,thick,scale=0.15mm]
\tikzstyle{cercle}=[circle,draw,thick,fill=black, scale=0.15mm]
\fill[color=black](8, 0)circle(0mm);
\fill[color=black](8, 6)circle(0mm);
\node[cercle] (C) at (16,15){};
\node[terminal] (T) at (0, 15){};
\draw[-](8, 0)--(8, 6);
\draw[-](8, 6)--(0, 14);
\draw[-](8, 6)--(16, 14);
\path (0,15) node[above]{$a$};
\path (16,15) node[above]{$b$};
\end{tikzpicture}}
\newcommand{\extthree}{\begin{tikzpicture}[ scale=0.05]
\tikzstyle{terminal}=[circle,draw,thick,scale=0.15mm]
\tikzstyle{cercle}=[circle,draw,thick,fill=black, scale=0.15mm]
\fill[color=black](8, 0)circle(0mm);
\fill[color=black](8, 6)circle(0mm);
\fill[color=black](4, 10)circle(0mm);
\node[cercle] (C) at (0,15){};
\node[cercle] (C) at (8,15){};
\node[terminal] (T) at (16, 15){};
\draw[-](8, 0)--(8, 6);
\draw[-](8, 6)--(0, 14);
\draw[-](4, 10)--(8, 14);
\draw[-](8, 6)--(16, 14);
\path(-2,15) node [above] {$-2$};
\path(8, 15) node[above] {$1$};
\path(16, 15) node[above] {$5$};
\end{tikzpicture}}
\newcommand{\extthreethree}{\begin{tikzpicture}[ scale=0.05]
\tikzstyle{terminal}=[circle,draw,thick,scale=0.15mm]
\tikzstyle{cercle}=[circle,draw,thick,fill=black, scale=0.15mm]
\fill[color=black](8, 0)circle(0mm);
\fill[color=black](8, 6)circle(0mm);
\fill[color=black](4, 10)circle(0mm);
\node[cercle] (C) at (0,15){};
\node[cercle] (C) at (8,15){};
\node[terminal] (T) at (16, 15){};
\draw[-](8, 0)--(8, 6);
\draw[-](8, 6)--(0, 14);
\draw[-](4, 10)--(8, 14);
\draw[-](8, 6)--(16, 14);
\path(-2,15) node [above] {$-3$};
\path(8, 15) node[above] {$0$};
\path(16, 15) node[above] {$4$};
\end{tikzpicture}}
\newcommand{\slfive}{\begin{tikzpicture}[ scale=0.05]
\tikzstyle{terminal}=[circle,draw,thick,scale=0.15mm]
\tikzstyle{cercle}=[circle,draw,thick,fill=black, scale=0.15mm]
\fill[color=black](9, -4)circle(0mm);
\fill[color=black](9, 2)circle(0mm);
\fill[color=black](4.5, 6)circle(0mm);
\fill[color=black](6, 8)circle(0mm);
\fill[color=black](15, 8)circle(1mm);
\node[cercle] (C) at (0,11.5){};
\node[cercle] (C) at (4,11.5){};
\node[cercle] (C) at (18,11.5){};
\node[terminal] (T) at (8, 11.5){};
\node[terminal] (T) at (12, 11.5){};
\draw[-](9, -4)--(9, 2);
\draw[-](9, 2)--(0, 11) ;
\draw[-](9, 2)--(18, 11);
\draw[-](4.5, 6)--(6, 8);
\draw[-](6, 8)--(4, 11) ;
\draw[-](6, 8)--(8, 11);
\draw[-](15, 8)--(12, 11);
\path (-2, 11.5)node [above] {$-2$};
\path (4,11.5) node [above] {$1$};
\path (8,11.5) node [above] {$5$};
\path(12,11.5)node [above] {$2$};
\path (18, 11.5) node[above] {$3$};
\end{tikzpicture}}

\title{Free Hom-groups, Hom-rings   and Semisimple modules}

 \author{ Imed Basdouri$^1$\footnote{Email : basdourimed@yahoo.fr}, Sami Chouaibi$^2$\footnote{Email : chouaibisami@yahoo.fr}, Abdenacer  Makhlouf$^3$\footnote{Email : abdenacer.makhlouf@uha.fr}, Esmael  Peyghan$^4$\footnote{Email : e-peyghan@araku.ac.ir}
  \\
 {\small $^1$  Faculty of Sciences, Department of Mathematics, Gafsa University, Tunisia.}\\
{\small $^2$  Faculty of Sciences, Department of Mathematics, Sfax University, Tunisia.}\\
{\small$^2$ Universit\'e de Haute Alsace, IRIMAS-d\'epartement de math\'ematiques, Mulhouse, France.}\\
{\small$^3$  Department of Mathematics, Faculty of Science, Arak University,
 	Arak,  Iran.}}
 
\begin{document} \maketitle

%

\begin{abstract}
The purpose of this paper is to introduce and study  a Hom-type generalization of rings. We provide  their   basic properties and and some key constructions. Furthermore, we consider modules over Hom-rings and  characterize the category of simple modules and simple Hom-rings. In addition, we extend some classical results and concepts of groups to Hom-groups. We construct free regular Hom-group using Super-Leaf weighted trees and discuss Normal Hom-subgroups,  abelianization of regular Hom-group, universality of tensor product of Hom-groups and   simple Hom-groups. 
\end{abstract}

\vspace{.5cm}



 \maketitle

\section*{Introduction}
Hom-type objects have been under intensive research in the last decade. The notion of Hom-Lie algebra was
introduced by Hartwig, Larsson, and Silvestrov in \cite{JDS} as part of a study of deformations of the Witt and the
Virasoro algebras \cite{NH},\cite{MPJ},\cite{TC}. Hom-associative algebras were introduced in \cite{AS}, where it is shown that the commutator bracket of a Hom-associative algebra gives rise to a Hom-Lie algebra and where a classification of Hom-Lie admissible algebras is established.
The first notion of Hom-group appeared first in \cite{CG}. Then it  was introduced as a non-associative analogue of a group in \cite{LMT}, where the authors first gave a
new construction of the universal enveloping algebra that is different from the one in \cite{DY}. This new construction
leads to a Hom-Hopf algebra structure on the universal enveloping algebra of a Hom-Lie algebra. One
can associate a Hom-group to any Hom-Lie algebra by considering group-like elements in its universal enveloping
algebra. Recently, M. Hassanzadeh developed representations and a (co)homology theory for Hom-groups in \cite{H1}.
He also proved Lagrange's theorem for finite Hom-groups in \cite{H2}. In \cite{JKS},  Hom-Lie groups and  relationship between Hom-Lie groups and
Hom-Lie algebras were explored. The author studied integrable Hom-Lie algebra describe a Hom-Lie group action on a smooth manifold. They  define a Hom-exponential (Hexp) map from the Hom-Lie algebra of a Hom-Lie group to the Hom-Lie group and discuss the universality of this Hexp map. We also describe a Hom-Lie group action on a smooth manifold. 

The purpose of this paper is to introduce Hom-rings which are twisted versions of the ordinary rings. We discuss
some of their properties and provide construction procedures using ordinary rings. Also, we  introduce the
notion of modules on  Hom-rings and characterize the category of simple modules and simple Hom-rings. Likewise, some developments of Hom-group theory are provide. In particular   free regular Hom-groups are constructed and normal Hom-subgroups explored.
The paper is organized as follows:
In Section 1, we recall some basic definitions and results concerning Hom-groups. In Section 2, we provide a construction of free regular Hom-group based on super-leaf weighted trees.  In Section 3, we introduce Normal Hom-subgroups and discuss their properties. In Section 4, we characterize simple Hom-groups and in Section 5 we deal with tensor products.In  Section 6, we introduce Hom-type versions of rings
and some related notions. Moreover, we study modules over
Hom-rings in Section 7.  Semisimple modules and Hom-rings are characterized in Section 8 and tensor product of modules discussed in the last section.

\section{Hom-groups}
We review some basics about  Hom-groups from \cite{LMT}, \cite{H1} and \cite{H2}, and provide  some new properties about hom-groups. Also, we recall  the notion of Hom-group homomorphism and  show a Hom-group structure on the set of all Hom-group homomorphisms between two abelian regular  Hom-groups. Finally, we discuss some classes of Hom-subgroups.
\subsection{Definitions and Examples}
\begin{df}\label{IM}
A Hom-group is a quadruplet $(G,\mu,e,\alpha)$ consisting of a set $G$ with a distinguished element $e$ of $G$, an operation $\mu:G\times G\rightarrow G$ and a set map $\a: G\rightarrow G$, such that the following  axioms are satisfied:
\begin{enumerate}
  \item The product map $\mu:G\times G\rightarrow G$ and the set map $\a: G\rightarrow G$  satisfy the Hom-associativity condition
  \begin{equation}\label{HOM1}
    \mu(\a(g),\mu(h,k))=\mu(\mu(g,h),\a(z)).
  \end{equation}
  For simplicity when there is no confusion we omit the multiplication sign $\mu.$
  \item The map $\a$ is multiplicative, i.e., $\a(\mu(g,h))=\mu(\a(g),\a(h))$.
  \item  The element $e$ is called unit and it satisfies the Hom-unitarity condition
  \begin{equation}\label{HOM2}
  \mu(g,e)=\mu(e,g)=\a(g)\;\;\quad \a(e)=e.
  \end{equation}
  \item the map $g\longmapsto g^{-1}$ is an antimorphism; $\mu(g,h)^{-1}=\mu(h^{-1},g^{-1}).$
  \item  For any $g\in G$ there exists a natural number $k$  satisfying the Hom-invertibility condition
  \begin{equation}\label{Hom3}
  \a^k(\mu(g,g^{-1}))=\a^k(\mu(g^{-1},g))=e.
  \end{equation}
  The smallest such $k$ is called the invertibility index of $g.$
  \end{enumerate}
If $\a$ is invertible,  condition (\ref{Hom3}) can be simplified to the following condition:\\
For every element $g\in G$, there exists an element $g^{-1}$ which  $\mu(g,g^{-1})=\mu(g^{-1},g)=e.$\\
A Hom-group, such as $\a$ is invertible, is called  regular Hom-group (\cite{H2},\cite{JKS}). If moreover the product $\mu$ is commutative, we say that $G$ is an abelian (commutative) Hom-group. In this case, the product $\mu$ will be denoted  by $+$ and the element $e$ by $0.$ Therefore the  multiplicativity axiom will change to $\a(g+h)=\a(g)+\a(h),~\forall g,~h\in G.$
\end{df}
\begin{rems}\label{calcul1}\item
\begin{enumerate}
\item From the Hom-unitarity and the Hom-associativity conditions,  we can show that the multiplicativity of the map $\a$ is intuitive;
$\a(\mu(g,h))=\mu(e,\mu(g,h))=\mu(\a(e),\mu(g,h))=\mu(\mu(e,g),\a(h))=\mu(\a(g),\a(h)$.
\item Since we have the antimorphism $g\longmapsto g^{-1}$ , therefore by the definition, inverse of any
element $g\in G$ is unique although different elements may have different invertibility indexes.
\item The inverse of the  element $e$ of a Hom-group is itself because $\a(\mu(e,e))=\a(e)=e.$
\item For any Hom-group $(G,\mu,e,\alpha)$, we have $\a(g^{-1})=\a(g)^{-1}$  and if the invertibility index of $g$ is $k$,  then the invertibility index of $\a(g)$ is $k-1.$
\end{enumerate}
\end{rems}
\begin{ex}
Let $(G,\mu,e)$ be a group and $\a:G\rightarrow G$ be a group homomorphism. We define a new product $\mu_{\a}:G\times G\rightarrow G$ given by
$$\mu_{\a}(g,h)=\a(\mu(g,h))=\mu(\a(g),\a(h)).$$
Then $(G,\mu_{\a},e,\a)$ is a Hom-group denoted by $G_{\a}$. It is said the twist group of $G$.
\end{ex}
\begin{ex}\label{produit}
Let $(G,\mu_G,e_G,\a_G)$ and $(H,\mu_H,e_H,\a_H)$ be two Hom-groups. Then $(G\times H,\mu,e,\a)$ is a Hom-group such as:
\begin{itemize}
  \item the identity element $e=(e_G,e_H),$
  \item the map $\a$ is defined by $\a(g,h)=(\a_G(g),\a_H(h)),\;\forall (g,h)\in G\times H,$
  \item the multiplication is given by $\mu((g_1,h_1),(g_2,h_2))=(\mu_G(g_1,g_2),\mu_H(h_1,h_2)),\;\forall (g_i,h_i)\in G\times H, \;i=1,~2.$
\end{itemize}
\end{ex}
We have this new example.
\begin{ex}
  Let $(\Z,+)$ be the additive group of integers. We have $End(\Z)\simeq \Z$, i.e., any endomorphism $\a\in End(\Z)$ is completely determined by the element $a(1)\in\Z$ and for every $q\in\Z$ there exists an endomorphism $\a\in End(\Z)$ such that $\a(1)=q$. In other words, the correspondence $\a\mapsto\a(1)$ is an isomorphism between $End(\Z)$ and $\Z.$ Then, for every $q\in\Z$, we can define a Hom-group denoted by $(\Z,+_q,0,\a)$ such as $+_q$ is defined by $n+_q m=q(n+m),~\forall n,~m\in\Z$ and checks the following condition
  \begin{equation}
    qn+_q( m+_q p)=(n+_qm)+_q qp=q^2(n+m+p),\;\forall n,~m,~p\in\Z.
  \end{equation}
  The  Hom-group $(\Z,+_q,0,\alpha )$ is called the $q$-additive group of integers.
\end{ex}

\subsection{Some properties}
\begin{lem}
Let $(G,\cdot,e,\a)$ be a Hom-group and $g\in G$. If $g\cdot g=\a(g)$ and $\a^2(g)=g$, then $g=e.$
\end{lem}
\begin{prf}
Let $h\in G$ such that $h\cdot g=e$. Then
\begin{align*}
  \a^2(g)=e\cdot\a(g)=(h\cdot g)\cdot\a(g)=\a(h)\cdot(g\cdot g)=\a(h)\cdot\a(g)=\a(h\cdot g)=\a(e)=e.
\end{align*}
From $\a^2(g)=g$, we can find out that $g=e.$
\end{prf}
\begin{pro}\label{R1}
Let $(G,\mu,e,\alpha)$ be a Hom-group and  $g$ be an invertible element   of index $n$. If $\a^{i}(g)h=\a^i(g)k,~\forall i\geq n$, then we have $\a^2(h)=\a^2(k)$.
\end{pro}
\begin{prf}
Using Definition \ref{IM}, we get
\begin{align*}
\a^2(h)=e\a(h)=\a^i(g^{-1}g)\a(h)=(\a^i(g^{-1})\a^i(g)))\a(h)=\a^{i+1}(g^{-1})(\a^i(g)h),\ \ \ \forall i\geq n.
\end{align*}
Applying $\a ^i(g)h=\a ^i(g)k$, the above equation reduces to the following:
\begin{align*}
\a^2(h)=\a^{i+1}(g^{-1})(\a^i(g)k)=\a^i(g^{-1}g)\a(k)=e\a(k)=\a^2(k).
\end{align*}
\end{prf}
\begin{cor}
If $(G,\mu,e,\alpha)$ is a regular Hom-group, $g$ is an invertible element and $gh=gk$, then $h=k$.
\end{cor}
\begin{lem}\label{RR1}
Let $(G,\mu,e,\alpha)$ be a Hom-group, $g$ be an invertible element of index $n$ and $l$ be an invertible element of index $m$. If $(\a^{i}(g)h)(k\a^j(l))=(\a^i(g)k)(h\a^j(l)),~\forall i\geq n,~j\geq m$, then $\a^3(hk)=\a^3(kh)$.
\end{lem}
\begin{prf}
Using Definition \ref{IM}, we get
  \begin{align*}
  \a^3(hk) &= \a^3(h)\a^3(k) \\
     &=
     \Big(\a^{i+2}(g^{-1})\Big(\a^{i+1}(g)\a(h)\Big)\Big)\Big(\Big(\a(k)\a^{j+1}(l)\Big)\a^{j+2}(l^{-1})\Big) \\
     &=
     \Big(\a^{i+2}(g^{-1})\Big(\a^{i+1}(g)\a(h)\Big)\Big)\a\Big((k\a^{j}(l))\a^{j+1}(l^{-1})\Big)  \\
     &=
     \a^{i+3}(g^{-1}) \Big(\Big(\a^{i+1}(g)\a(h)\Big)\Big((k\a^{j}(l))\a^{j+1}(l^{-1})\Big)\Big) \\
     &=
      \a^{i+3}(g^{-1})\Big(\Big((\a^{i}(g)h)(k\a^{j}(l))\Big)\a^{j+2}(l^{-1})\Big) \\
     &=
     \a^{i+3}(g^{-1})\Big(\Big((\a^{i}(g)k)(h\a^{j}(l))\Big)\a^{j+2}(l^{-1})\Big)  \\
     &=
     \a^{i+3}(g^{-1}) \Big(\Big(\a^{i+1}(g)\a(k)\Big)\Big((h\a^{j}(l))\a^{j+1}(l^{-1})\Big)\Big) \\
     &=
    \a^{i+3}(g^{-1})\Big(\Big(\a^{i+1}(g)\a(k)\Big)\a^2(h)\Big) \\
     &=
     \a^{i+3}(g^{-1})\Big(\a^{i+2}(g)\a(kh)\Big)\Big) \\
     &=
     \a^{i+2}(g^{-1}g)\a^2(kh)\\
     &=
     \a^3(kh).
  \end{align*}
\end{prf}
\begin{cor}
Let $(G,\mu,e,\alpha)$ be a regular Hom-group and $g$, $l$ be invertible elements. If $(gh)(kl)=(gk)(hl)$, then $hk=kh$.
\end{cor}
\begin{pro}\label{RRR1}
Let $(G,\mu,e,\alpha)$ be a Hom-group, $g$ be an invertible element of index $n$ and $l$ be an invertible element of index $m$.
If $(\a^{i}(g^{-1})\a^{j}(h^{-1}))(\a^{i}(g)\a^{j}(h))=e,~\forall i\geq n,~j\geq m$, then $\a^{i+5}(g)\a^{j+5}(h)=\a^{j+5}(h)\a^{i+5}(g).$
\end{pro}
\begin{prf}
We have $e=\a^{i}(g^{-1}g)\a^{j}(h^{-1}h)$. Lemma \ref{RR1} gives us $\a^3(\a^i(g)\a^j(h^{-1}))=\a^3(\a^j(h^{-1})\a^i(g)).$\\
By multiplying this last equation on the  left by $\a^{j++4}(h)$ and using the hom-associativity of $\a$, we get
$$\a^{i+5}(g)=\a^{j+4}(h)\Big(\a^{i+3}(g)\a^{j+3}(h^{-1})\Big)=\Big(\a^{j+3}(h)\a^{i+3}(g)\Big)\a^{j+3}(h^{-1}).$$
Multiply this equation on the  right by $\a^{j++5}(h)$, we find $$\a^{i+5}(g)\a^{j+5}(h)=\a^{j+5}(h)\a^{i+5}(g).$$
\end{prf}
\begin{cor}
Let $(G,\mu,e,\alpha)$ be a regular Hom-group and $g$, $l$ be  invertible elements of $G$. If  $(g^{-1}h^{-1})(gh)=e$, then $gh=hg$.
\end{cor}
\begin{pro}\label{lemmeRS}
Let $(G,\cdot,e,\alpha)$ be a regular abelian Hom-group. For all $g,~h,~k$ and $l$ in $G$, we have
\begin{equation}\label{lemmers}
  (g\cdot h)\cdot(k\cdot l)=(g\cdot k)\cdot(h\cdot l).
\end{equation}
\end{pro}
\begin{prf}
For all $g,~h,~k$ and $l$ in $G$, we have
\begin{align*}
  (g\cdot h)\cdot(k\cdot l) &=\a(g)\cdot(h\cdot\a^{-1}(k\cdot l))=
   \a(g)\cdot(\a^{-1}(h\cdot k)\cdot l)=
   \a(g)\cdot(\a^{-1}(k\cdot h)\cdot l)\\
   &=
   \a(g)\cdot(k\cdot\a^{-1}(h\cdot l))=
   (g\cdot k)\cdot(h\cdot l).
\end{align*}
\end{prf}
\subsection{Hom-groups Homomorphisms}
\begin{df}
Let $(G,\mu_G,e_G,\alpha_G)$ and $(H,\mu_H,e_H,\alpha_H)$ be two Hom-groups.
A homomorphism of Hom-groups is a map $f: G\rightarrow H$ such that $$f(\mu_G(g,h))=\mu_H(f(g),f(h)),~\forall g,~h \in G\;\; \text{ and } \;\; f\circ\a_G=\a_H\circ f.$$
It is said to be a weak homomorphism if $f\circ\a_G=\a_H\circ f$.
\end{df}
\begin{rem}
Let $f: G\rightarrow H$ be a homomorphism of Hom-groups such that $f(e_G)=e_H$  and $g$ be an element in $G$. Then, we have
\begin{align*}
  f(\a_G(g)) =f(\mu_G(e_G,g))=
   \mu_H(f(e_G),f(g))=
   \mu_H(e_H,f(g))=
   \a_H(f(g)).
\end{align*}
So any homomorphism of Hom-groups $f:(G,\mu_G,e_G,\alpha_G)\rightarrow (H,\mu_H,e_H,\alpha_H)$ such that $f(e_G)=e_H$ is a weak homomorphism.
\end{rem}
\begin{df}
Two Hom-groups $(G,\mu_G,e_G,\alpha_G)$ and $(H,\mu_H,e_H,\alpha_H)$ are called isomorphic if there exists a bijective
homomorphism of Hom-groups $f:(G,\mu_G,e_G,\alpha_G)\rightarrow (H,\mu_H,e_H,\alpha_H)$.
\end{df}
\begin{pro}
Let $(G,\cdot,e,\alpha)$ be a regular Hom-group and define a map $f:G\rightarrow G$ by $f(g)=g\cdot g$ for each $g\in G$.
Then $G$ is an abelian Hom-group if and only if the map $f$ is a Hom-group homomorphism.
\end{pro}
\begin{prf}
 Suppose that $G$ is an abelian Hom-group. Let $g,~h$ be arbitrary elements in $G$. Then using Lemma \ref{lemmeRS}, we have
  \begin{align*}
    f(g\cdot h)=(g\cdot h)\cdot (g\cdot h)=(g\cdot g)\cdot (h\cdot h)=f(g)\cdot f(h).
  \end{align*}
 Hence $f$ is a Hom-group homomorphism from $G$ to $G$. Conversely, suppose that $f:G\rightarrow G$ by $f(g)=g\cdot g$ is a Hom-group homomorphism. Let $g,~h$ be arbitrary elements in $G$. We prove that $g\cdot h=h\cdot g$. Since $f$ is a Hom-group homomorphism, we have $f(g\cdot h)=f(g)\cdot f(h),$ i.e., $(g\cdot h)\cdot (g\cdot h)=(g\cdot g)\cdot (h\cdot h)$.
 By Hom-associativity of the multiplication of $G$ , the last equation is equivalent to $\a(g)\cdot(h\cdot\a^{-1}(g\cdot h))=\a(g)\cdot(g\cdot\a^{-1}(h\cdot h)).$ Multiplying this by $\a^2(g^{-1})$ on the left, we obtain $\a^2(h)\cdot\a(g\cdot h)= \a^2(g)\cdot\a(h\cdot h)$. By Hom-associativity we have $\a(h\cdot g)\cdot\a^2(h)= \a(g\cdot h)\cdot\a^2(h)$. Multiplying this by $\a^3(h^{-1})$ on the right, we find $\a^3(h\cdot g)=\a^3(g\cdot h).$ Since $\a$ is bijective, we have $g\cdot h=h\cdot g.$ So $G$ is abelian.
\end{prf}
\begin{pro}
Let  $(G,+_G,e_G,\alpha_G)$ and $(H,+_H,e_H,\alpha_H)$ be abelian regular Hom-groups.\\  Then $(Hom(G, H),+,0,\alpha)$ is  an abelian Hom-group, where $Hom(G, H)$ is the set  of all Hom-group homomorphisms from $G$ to $H$,  $f+g$ is defined by
\begin{equation}\label{homomorgp}
  (f+g)(x)=f(x)+_Hg(x),\;\;\mbox{for all}\;\;x\in G,
\end{equation}
and $\a:Hom(G, H)\rightarrow Hom(G, H)$ is defined by $\a(f)(x)=\a_H(f(x))$, for all $x\in G.$
\end{pro}
\begin{prf}
According to Lemma \ref{lemmeRS} and from the commutativity of $H$, we can prove that $f+g$ is again a Hom-group homomorphism. The other properties for $Hom(G, H)$ to be an abelian Hom-group are consequences of the fact that $(H,+_H,e_H,\alpha_H)$ is an abelian Hom-group.
\end{prf}
\subsection{Hom-subgroups}
\begin{df}
Let $(G,\mu,e,\alpha)$ be a Hom-group. A nonempty subset $H$ of the Hom-group $G$ is a Hom-subgroup of $G$ if $(H,\mu,e,\alpha)$ is a Hom-group. We use the notation $H\leq G$ to indicate that $H$ is a Hom-subgroup of $G.$
\end{df}
If $H$ is a Hom-subgroup, we see that $e$ the identity for $G$ is also
the identity for $H$. Consequently the following theorem is obvious:
\begin{thm}
A subset $H$ of the Hom-group $(G,\cdot,e,\alpha)$ is a Hom-subgroup of $G$ if and only if
\begin{enumerate}
  \item $e\in H,$
  \item $h\in H \Rightarrow h^{-1}\in H,$
  \item $h_1,~h_2\in H \Rightarrow h_1\cdot h_2\in H$.
\end{enumerate}
\end{thm}
\begin{rems}
\begin{enumerate}
\item From the third item of the previous theorem, one notes that if $H$ is a Hom-subgroup of $G$ then $\a(h)=e\cdot h\in H$ for all $h\in H$. Therefore $\a(H)\subset H.$
\item Let $(G,\mu,e,\alpha)$ be a Hom-group and $H$ be a subset of $G$. A Hom-group structure on H is defined by the inclusion map.
\end{enumerate}
\end{rems}
\begin{ex}\label{exex}
Let $G$ be a group and $\a:G\rightarrow G$ be a group homomorphism. If $H$ is a subgroup of $G$ which $\a(H)\subset H$, then $H_{\a}$ is a Hom-subgroup of $G_{\a}.$
\end{ex}
\begin{ex}
Let $(G,\mu_G,e_G,\alpha_G)$ and $(H,\mu_H,e_H,\alpha_H)$. If $f:G\rightarrow H$ is a Hom-group homomorphism such that $f(e_G)=e_H$, then the inverse image of each Hom-subgroup of $H$  is a Hom-subgroup of $G.$
\end{ex}
\begin{df}\cite{H1}
Let $(G,\mu,e,\alpha)$ be a regular Hom-group. The center $Z(G)$ of the Hom-group $(G,\mu,e,\alpha)$ is the set of all $g\in G$ where $gx=xg$ for all $x\in G.$
\end{df}
\begin{pro}\cite{H1}
$(Z(G),\mu,e,\alpha)$ is a Hom-subgroup  of $(G,\mu,e,\alpha)$.
\end{pro}
\begin{lem}\label{Z}
Let $(G,\mu,e,\alpha)$ be a regular Hom-group. Then, $\forall g\in Z(G),$ we have $\a^k(g)\in Z(G),~\forall k\in\N.$
\end{lem}
\begin{prf}
Let $g\in Z(G)$.
\begin{enumerate}
  \item For $k=1$, we have  $\a(g)x=(eg)x=\a(e)(g\a^{-1}(x))=e(\a^{-1}(x)g)=(e\a^{-1}(x))\a(g)=x\a(g).$ Then $\a(g)\in Z(G)$.
  \item Suppose that for any $k\in\N$ we have $\a^k(g)\in Z(G)$ and show that $\a^{k+1}(g)$ belongs to $Z(G)$;
  $$\a^{k+1}(g)x=(e\a^k(g))x=\a(e)(\a^k(g)\a^{-1}(x))=e(\a^{-1}(x)\a^k(g))=(e\a^{-1}(x))\a^{k+1}(g)=x\a^{k+1}(g).$$
\end{enumerate}
\end{prf}
Let $(G,\mu,e)$ be a group and $(G,\mu_{\a},e,\alpha)$ the associated twist group of $(G,\mu,e)$. Consider the set
$$Z_{\a}(G)=\{g\in G/ \a(gx)=\a(xg),~\;\forall x\in G\}.$$ Then, we have the following proposition
\begin{pro}
\begin{enumerate}
If $Z(G)$ is  the center of the group $(G,\mu,e)$, so we have
  \item $Z(G)\subset \a^{-1}(\a(Z(G)))\subset Z_{\a}(G).$
  \item $Z(G)=Z_{\a}(G)$ if $\a$ is injective.
  \item If $\a(Z(G))\subset Z(G)$ and $\a$ is injective, then $Z_{\a}(G)$ is a Hom-subgroup of $(G,\mu_{\a},e,\alpha).$
\end{enumerate}
\end{pro}
\begin{prf}
\begin{enumerate}
  \item Let $g\in \a^{-1}(\a(Z(G)))$. Then, it exists $h\in Z(G)$ such as $\a(g)=\a(h).$ Let us show that $g$ is an element in $Z_{\a}(G)$, i.e., $\a(gx)=\a(xg),~\;\forall x\in G.$
\begin{align*}
 \a(gx )=\a(g)\a(x)=\a(h)\a(x)=\a(hx)=\a(xh)=\a(x)\a(h)=\a(x)\a(g)=\a(xg).
\end{align*}
This shows that $g\in Z_{\a}(G).$
  \item Let $g\in Z_{\a}(G)$. So, $ \a(gx )=\a(xg)$ for all $x\in G$. Since $\a$ is injective, then $gx=xg,~\forall x\in G$. Therefore, $ Z_{\a}(G)\subset Z(G).$
  \item According to  Example \ref{exex} and the previous property, we can show that $Z_{\a}(G)$ is a Hom-subgroup of $(G,\mu_{\a},e,\alpha)$.
\end{enumerate}
\end{prf}
\section{Free regular Hom-group}
\subsection{Free regular Hom-group with $1-$generator}
A planar tree is an oriented graph drawn on a plane with only one root. It is called binary when any vertex is trivalent, i.e.,
one root and two leaves. Usually we draw the root at the bottom of the tree and the leaves are drawn at
the top of it:
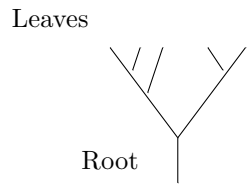
\begin{align*}
\begin{tikzpicture}[ scale=0.1]
\fill[color=black](9, -4)circle(0mm);
\fill[color=black](9, 2)circle(0mm);
\fill[color=black](0, 14)circle(0mm);
\fill[color=black](18, 14)circle(0mm);
\fill[color=black](3, 11)circle(0mm);
\fill[color=black](4, 14)circle(0mm);
\fill[color=black](15, 11)circle(0mm);
\fill[color=black](13, 14)circle(0mm);
\fill[color=black](5, 8)circle(0mm);
\fill[color=black](7, 14)circle(0mm);
\draw[-](9, -4)--(9, 2);
\draw[-](9, 2)--(0, 14);
\draw[-](9, 2)--(18, 14);
\draw[-](3, 11)--(4, 14);
\draw[-](15, 11)--(13, 14);
\draw[-](5, 8)--(7, 14);
\path(0,-1) node {Root};
\path(-8,18) node {Leaves};
\end{tikzpicture}
\end{align*}
For any natural number $ n\geq1$, let $T_n$ denote the set of planar binary trees with $n$ leaves and one root. The first $T_n$ are depicted below.\\
\[
   T_1 = \left\lbrace \tone \, \right\rbrace,~
   T_2 = \left\lbrace \ttwo \right\rbrace,~
   T_3 = \left\lbrace \tthreeone,\, \tthreetwo \right\rbrace,~
   T_4 = \left\lbrace \tfourone,\, \tfourtwo,\, \tfourthree,\, \tfourfour,\, \tfourfive \right\rbrace.
   \]

An element in $T_n$ will be called an $n$-tree. When necessary we label the leaves of an $n$-tree by
$1,~2,~3,~ \cdot\cdot\cdot,~ n$ from left to right.\\
Let $\psi\in T_n$ and $\varphi\in T_m$ be a pair of trees, the $(n + m)$-tree $\psi\vee\varphi$, called the grafting of $\psi$ and $\varphi$, is
obtained by joining the roots of $\psi$ and $\varphi$ to create a new root. For instance,
$$\psi\vee\varphi=\psiwedgephi.$$
\begin{rem}
Note that grafting is a nonassociative and non-commutative operation. For any tree $\varphi\in T_n$, there are unique integers $p$ and $q$ with $p+q=n$ and trees $\varphi_1\in T_p$ and $\varphi_2\in T_q$ such that $\varphi=\varphi_1\vee\varphi_2$. The tree $\varphi_1$ is called the left subtree (L-subtree) of $\varphi$ and $\varphi_2$ is the right subtree (R-subtree) of $\varphi$. So any tree $\varphi$ will be represented as follows
$$\varphi=\ntree.$$
\end{rem}
\begin{df}
Let $\varphi=\varphi_1\vee\varphi_2$ be a tree. Two successive leaves of $\varphi$  are called
\begin{enumerate}
  \item left (resp. right) attached leaves (LAL resp. RAL) if they have the same node in the left (resp. right) subtree of $\varphi$:
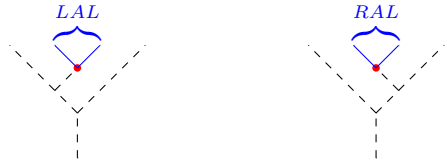
   $$\LAL\quad\quad\quad\quad\quad\quad\RAL$$
  \item left (resp. right) disjoint leaves (LDL resp. RDL) if they don't have the same node in the left (resp. right) subtree of $\varphi$: $$\LDL\quad\quad\quad\quad\quad\quad\RDL$$
  \item disjoint leaves (DL) if the first leaf is the last leaf of the left subtree of $\varphi$ and the other is the first leaf of the right subtree of $\varphi$:
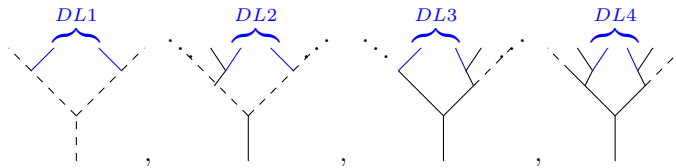
$$\DLon,\DLtw,\DLtr,\DLfr$$
\end{enumerate}
\end{df}
\begin{df}
A bicolored Leaf $n$-tree is a $n$-tree such that each leaf is colored with black (b) or white (w). The set of bicolored Leaf $n$-tree will be denoted by $D_n$.
\end{df}
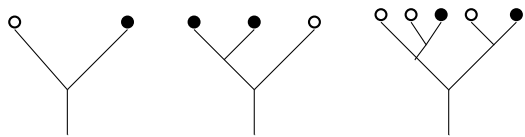
\begin{exs}
$$\begin{matrix}
  \exttwoD\quad\quad & \extthreeD\quad\quad & \extfiveD
\end{matrix}$$
\end{exs}
\begin{df}
A Super-Leaf weighted $n$-tree is a pair $(\varphi,a)$ where:
\begin{itemize}
  \item $\varphi$ is an element of $D_n$ ,
  \item $a$ is an $n-$tuplet $(a_1,a_2,\cdot\cdot\cdot,a_n)$  such that, for all $1\leq i\leq n,~a_i\in\Z.$
\end{itemize}
We call the bicolored Leaf tree $\varphi$ the underlying tree of the Super-Leaf weighted $n$-tree $(\varphi,a)$. For all $1\leq i\leq n$, the integer $a_i$ is said to be the weight of the leaf $i$.
\end{df}
We will indeed barely use the notation $(\varphi,a)$  at all, and find more convenient to picture a Super-leaf weighted $n$-tree $(\varphi,a_1,a_2,\cdot\cdot\cdot,a_n)$ by drawing the tree $\varphi$ and putting the weight $a_i$ next to each leaf. For example, here are a Super-leaf weighted $2$-tree and a Super-leaf weighted $3$-tree:
$$\exttwo\quad\quad\quad\quad\extthree$$
The grafting operation extends to leaf weighted $n$-trees. For example:
$$\extthree\vee\exttwo=\slfive$$
For $n\geq1$, we let $B_n$ denote the set of all Super-Leaf weighted $n$-trees. Let $B$ denote the union over $n\in\N$ of the sets $B_n$ together with an element that we call the unit and denote by $\mathbf{1}$. Note that the element $\mathbf{1}$ is different from the Super-Leaf weighted $1$-tree
$$
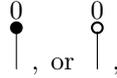
\begin{tikzpicture}[ scale=0.05]
\tikzstyle{terminal}=[circle,draw,thick,scale=0.15mm]
\tikzstyle{cercle}=[circle,draw,thick,fill=black, scale=0.15mm]
\fill[color=black](1, 0)circle(0mm);
\node[cercle] (C) at (1,11){};
\draw[-](1, 0)--(1, 10);
\path (1,11)node[above]{$0$};
\end{tikzpicture},\;\mbox{or}\; \begin{tikzpicture}[ scale=0.05]
\tikzstyle{terminal}=[circle,draw,thick,scale=0.15mm]
\tikzstyle{cercle}=[circle,draw,thick,fill=black, scale=0.15mm]
\fill[color=black](1, 0)circle(0mm);
\node[terminal] (T) at (1,11){};
\draw[-](1, 0)--(1, 10);
\path (1,11)node[above]{$0$};
\end{tikzpicture},$$
and which checks the following property
$$\mathbf{1}=
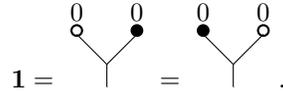
\begin{tikzpicture}[ scale=0.05]
\tikzstyle{terminal}=[circle,draw,thick,scale=0.15mm]
\tikzstyle{cercle}=[circle,draw,thick,fill=black, scale=0.15mm]
\fill[color=black](8, 0)circle(0mm);
\fill[color=black](8, 6)circle(0mm);
\node[cercle] (C) at (16,15){};
\node[terminal] (T) at (0, 15){};
\draw[-](8, 0)--(8, 6);
\draw[-](8, 6)--(0, 14);
\draw[-](8, 6)--(16, 14);
\path (0,15) node[above]{$0$};
\path (16,15) node[above]{$0$};
\end{tikzpicture}=\begin{tikzpicture}[ scale=0.05]
\tikzstyle{terminal}=[circle,draw,thick,scale=0.15mm]
\tikzstyle{cercle}=[circle,draw,thick,fill=black, scale=0.15mm]
\fill[color=black](8, 0)circle(0mm);
\fill[color=black](8, 6)circle(0mm);
\node[cercle] (C) at (0,15){};
\node[terminal] (T) at (16, 15){};
\draw[-](8, 0)--(8, 6);
\draw[-](8, 6)--(0, 14);
\draw[-](8, 6)--(16, 14);
\path (0,15) node[above]{$0$};
\path (16,15) node[above]{$0$};
\end{tikzpicture}.$$
We then consider the free vector space $\mathbb{T}$ generated by the set $B$.\\
We define on $\mathbb{T}$ two natural linear maps:
\begin{enumerate}
  \item $\a:\mathbb{T}\rightarrow \mathbb{T} $ sending $\mathbf{1}$ to $\mathbf{1}$ and sending a Super-Leaf weighted $n$-tree $(\varphi,a_1,a_2,\cdot\cdot\cdot,a_n)$ to the Super-Leaf weighted $n$-trees $(\varphi,a_1+1,a_2+2,\cdot\cdot\cdot,a_n+1).$
    \item a product $\vee$ that, for any pair of Super-Leaf weighted trees, is just the grafting of these trees and such
that for any Super-Leaf weighted $=\varphi=(\varphi,a)$
\begin{equation}\label{RFD2}
  \mathbf{1}\vee\varphi=\varphi\vee\mathbf{1}=\a(\varphi)
\end{equation}
\end{enumerate}
For examples
\begin{enumerate}
  \item $\a\Big(\exttwotwo\Big)=\exttwo$
  \item $\a\Big(\extthreethree\Big)=\extthree.$
\end{enumerate}
\begin{rem}
  Since $\a(\mathbf{1})=\mathbf{1}$, we can notice that, for all $a\in\Z$;
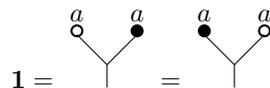
  $$\mathbf{1}=\begin{tikzpicture}[ scale=0.05]
\tikzstyle{terminal}=[circle,draw,thick,scale=0.15mm]
\tikzstyle{cercle}=[circle,draw,thick,fill=black, scale=0.15mm]
\fill[color=black](8, 0)circle(0mm);
\fill[color=black](8, 6)circle(0mm);
\node[cercle] (C) at (16,15){};
\node[terminal] (T) at (0, 15){};
\draw[-](8, 0)--(8, 6);
\draw[-](8, 6)--(0, 14);
\draw[-](8, 6)--(16, 14);
\path (0,15) node[above]{$a$};
\path (16,15) node[above]{$a$};
\end{tikzpicture}=\begin{tikzpicture}[ scale=0.05]
\tikzstyle{terminal}=[circle,draw,thick,scale=0.15mm]
\tikzstyle{cercle}=[circle,draw,thick,fill=black, scale=0.15mm]
\fill[color=black](8, 0)circle(0mm);
\fill[color=black](8, 6)circle(0mm);
\node[cercle] (C) at (0,15){};
\node[terminal] (T) at (16, 15){};
\draw[-](8, 0)--(8, 6);
\draw[-](8, 6)--(0, 14);
\draw[-](8, 6)--(16, 14);
\path (0,15) node[above]{$a$};
\path (16,15) node[above]{$a$};
\end{tikzpicture}$$
\end{rem}
The proof of the following lemma is obvious.
\begin{lem}
The map $\a$ defined above is a morphism for the grafting of trees, that is
$$\a(\varphi\vee\psi)=\a(\varphi)\vee\a(\psi),\;\mbox{for any}\;\varphi,\psi\in\mathbb{T}.$$
In addition, $\a$ is bijective and $\a^{-1}$ is defined, for all $(\varphi,a_1,a_2,\cdot\cdot\cdot,a_n)\in\mathbb{T}$, by\\
$$\a^{-1}(\varphi,a_1,a_2,\cdot\cdot\cdot,a_n)=(\varphi,a_1-1,a_2-1,\cdot\cdot\cdot,a_n-1).$$
\end{lem}
\begin{df}
Two leaves of a Super-Leaf weighted $n$-tree are said inverse type leaves if one is colored by b then the other is colored by w.
\end{df}
We now define a reduction process which allows one to obtain a reduced Super-Leaf weighted $\Omega$-type $n$-trees from an arbitrary Leaf weighted $\Omega$-type $n$-tree.
\begin{df}\label{reduction}\item
Let $\varphi $ be a Super-Leaf weighted $n$-tree. An elementary reduction of $\varphi$ is defined as follows
\begin{enumerate}
  \item If the leaf $i$ and the leaf $i+1$ are two left attached  and inverse type leaves of the same weight, then an elementary reduction of $\varphi$ is obtained by replacing the leaf $i$ and the leaf $i+1$ by $\mathbf{1}.$
  \item If the leaf $i$ and the leaf $i+1$ are two right attached  and inverse type leaves of the same weight, then an elementary reduction of $\varphi$ is obtained by replacing the leaf $i$ and the leaf $i+1$ by $\mathbf{1}.$
  \item If the leaf $i$ and the leaf $i+1$ are two left disjoint  and inverse type leaves of  weights respectively $b$ and $b+1$, then an elementary reduction of $\varphi$ is obtained by replacing the leaf $i$ and the leaf $i+1$ by $\mathbf{1}.$
  \item If the leaf $i$ and the leaf $i+1$ are two right disjoint  and inverse type leaves of  weights respectively $b+1$ and $b$, then an elementary reduction of $\varphi$ is obtained by replacing the leaf $i$ and the leaf $i+1$ by $\mathbf{1}.$
  \item If the leaf $i$ and the leaf $i+1$ are two DL1  and inverse type leaves of the same weight, then an elementary reduction of $\varphi$ is obtained by replacing the leaf $i$ and the leaf $i+1$ by $\mathbf{1}.$
  \item If the leaf $i$ and the leaf $i+1$ are two DL2  and inverse type leaves of weights respectively $b$ and $b+1$, then an elementary reduction of $\varphi$ is obtained by replacing the leaf $i$ and the leaf $i+1$ by $\mathbf{1}.$
  \item If the leaf $i$ and the leaf $i+1$ are two DL3  and inverse type leaves of weights respectively $b+1$ and $b$, then an elementary reduction of $\varphi$ is obtained by replacing the leaf $i$ and the leaf $i+1$ by $\mathbf{1}.$
  \item If the leaf $i$ and the leaf $i+1$ are two DL4  and inverse type leaves of of the same weight, then an elementary reduction of $\varphi$ is obtained by replacing the leaf $i$ and the leaf $i+1$ by $\mathbf{1}.$
\end{enumerate}
\end{df}
\begin{exs}
\begin{enumerate}
  \item Suppose $\varphi=(\varphi_1\vee\extredtwon)\vee\varphi_2$ for some Super-Leaf weighted trees $\varphi_1,~\varphi_2$ and $a\in\Z$. Then the elementary reduction of $\varphi$ results in the Super-Leaf weighted tree $\a(\varphi_1)\vee\varphi_2.$
  \item If $\varphi=\varphi_1\vee(\varphi_2\vee\extredtwon)$ for some Super-Leaf weighted trees $\varphi_1,~\varphi_2$ and $a\in\Z$. Then the elementary reduction of $\varphi$ results in the Super-Leaf weighted tree $\varphi_1\vee\a(\varphi_2).$
  \item Suppose $\varphi=((\varphi_1\vee\redoneone)\vee\redtwtw)\vee\varphi_2$ for some Super-Leaf weighted trees $\varphi_1,~\varphi_2$ and $a\in\Z$. Then the elementary reduction of $\varphi$ results in the Super-Leaf weighted tree $\a^2(\varphi_1)\vee\varphi_2.$
  \item If $\varphi= \varphi_1\vee(\redonetwone\vee(\redonetw\vee\varphi_2))$ for some Super-Leaf weighted trees $\varphi_1,~\varphi_2$ and $a\in\Z$. Then the elementary reduction of $\varphi$ results in the Super-Leaf weighted tree $\varphi_1\vee\a^2(\varphi_2).$
  \item Suppose $\varphi=(\varphi_1\vee\redoneone)\vee(\redonetw\vee\varphi_2)$ for some Super-Leaf weighted trees $\varphi_1,~\varphi_2$ and $a\in\Z$. Then the elementary reduction of $\varphi$ results in the Super-Leaf weighted tree $\a(\varphi_1)\vee\a(\varphi_2).$
  \item Suppose $\varphi=(\varphi_1\vee\extredtwtw)\vee(\redtwtw\vee\varphi_2)$ for some Super-Leaf weighted trees $\varphi_1,~\varphi_2$ and $a,b\in\Z$. Then the elementary reduction of $\varphi$ results in the Super-Leaf weighted tree $(\varphi_1\vee\begin{tikzpicture}[ scale=0.05]
\tikzstyle{terminal}=[circle,draw,thick,scale=0.15mm]
\tikzstyle{cercle}=[circle,draw,thick,fill=black, scale=0.15mm]
\fill[color=black](1, 0)circle(0mm);
\node[terminal] (T) at (1,11){};
\draw[-](1, 0)--(1, 10);
\path (1,11)node[above]{$b+1$};
\end{tikzpicture})\vee\a(\varphi_2).$
  \item If $\varphi=(\varphi_1\vee\redonetwone)\vee(\extredtwtwthr\vee\varphi_1)$ for some Super-Leaf weighted trees $\varphi_1,~\varphi_2$ and $a,b\in\Z$. Then the elementary reduction of $\varphi$ results in the Super-Leaf weighted tree $\a(\varphi_1)\vee (\begin{tikzpicture}[scale=0.05]
\tikzstyle{terminal}=[circle,draw,thick,scale=0.15mm]
\tikzstyle{cercle}=[circle,draw,thick,fill=black, scale=0.15mm]
\fill[color=black](1, 0)circle(0mm);
\node[cercle] (C) at (1,11){};
\draw[-](1, 0)--(1, 10);
\path (1,11)node[above]{$b+1$};
\end{tikzpicture}\vee\varphi_2)$.
  \item Suppose $\varphi=(\varphi_1\vee\begin{tikzpicture}[ scale=0.03]
\tikzstyle{terminal}=[circle,draw,thick,scale=0.15mm]
\tikzstyle{cercle}=[circle,draw,thick,fill=black, scale=0.15mm]
\fill[color=black](8, 0)circle(0mm);
\fill[color=black](8, 6)circle(0mm);
\node[cercle] (C) at (16,15){};
\node[terminal] (T) at (0, 15){};
\draw[-](8, 0)--(8, 6);
\draw[-](8, 6)--(0, 14);
\draw[-](8, 6)--(16, 14);
\path (0,15) node[above]{$b$};
\path (16,15) node[above]{$a$};
\end{tikzpicture})\vee(\begin{tikzpicture}[ scale=0.03]
\tikzstyle{terminal}=[circle,draw,thick,scale=0.15mm]
\tikzstyle{cercle}=[circle,draw,thick,fill=black, scale=0.15mm]
\fill[color=black](8, 0)circle(0mm);
\fill[color=black](8, 6)circle(0mm);
\node[terminal] (T) at (16,15){};
\node[terminal] (T) at (0, 15){};
\draw[-](8, 0)--(8, 6);
\draw[-](8, 6)--(0, 14);
\draw[-](8, 6)--(16, 14);
\path (0,15) node[above]{$a$};
\path (16,15) node[above]{$c$};
\end{tikzpicture}\vee\varphi_2)$ for some Super-Leaf weighted trees $\varphi_1,~\varphi_2$ and $a,b,~c\in\Z$. Then the elementary reduction of $\varphi$ results in the Super-Leaf weighted tree $(\varphi_1\vee\begin{tikzpicture}[ scale=0.05]
\tikzstyle{terminal}=[circle,draw,thick,scale=0.15mm]
\tikzstyle{cercle}=[circle,draw,thick,fill=black, scale=0.15mm]
\fill[color=black](1, 0)circle(0mm);
\node[terminal] (T) at (1,11){};
\draw[-](1, 0)--(1, 10);
\path (1,11)node[above]{$b+1$};
\end{tikzpicture})(\begin{tikzpicture}[ scale=0.05]
\tikzstyle{terminal}=[circle,draw,thick,scale=0.15mm]
\tikzstyle{cercle}=[circle,draw,thick,fill=black, scale=0.15mm]
\fill[color=black](1, 0)circle(0mm);
\node[terminal] (T) at (1,11){};
\draw[-](1, 0)--(1, 10);
\path (1,11)node[above]{$c+1$};
\end{tikzpicture}\vee\varphi_2).$
\end{enumerate}
\end{exs}
\begin{df}
Let $\varphi$ be a Super-Leaf weighted tree. A reduction of $\varphi$ (or a reduction process starting at $\varphi$) consists of consequent applications of elementary reductions starting at $\varphi$ and ending at a reduced Super-Leaf weighted tree:
\begin{align*}
  \varphi\rightarrow\varphi_1\rightarrow\varphi_2\rightarrow\cdots\rightarrow\varphi_n.
\end{align*}
The Super-Leaf weighted tree $\varphi_n$ is termed the reduced Super-Leaf weighted tree of $\varphi$ and it is denoted $\overline{\varphi}.$
\end{df}
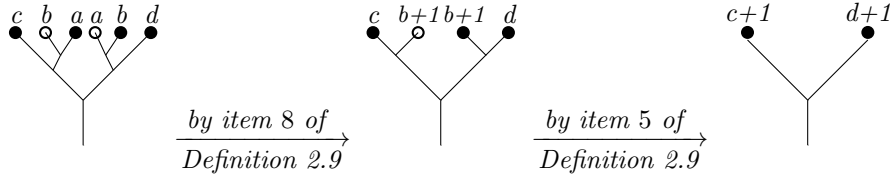
\begin{ex}\item
\begin{align*}
\begin{tikzpicture}[ scale=0.1]
\tikzstyle{terminal}=[circle,draw,thick,scale=0.15mm]
\tikzstyle{cercle}=[circle,draw,thick,fill=black, scale=0.15mm]
\fill[color=black](9, -4)circle(0mm);
\fill[color=black](9, 2)circle(0mm);
\fill[color=black](6, 8)circle(0mm);
\fill[color=black](4, 11)circle(0mm);
\fill[color=black](5,6)circle(0mm);
\fill[color=black](13, 6)circle(0mm);
\fill[color=black](12, 8)circle(0mm);
\node[cercle] (C) at (0,11){};
\node[terminal] (T) at (4,11){};
\node[cercle] (C) at (8,11){};
\node[terminal] (T) at (10.6,11){};
\node[cercle] (C) at (14,11){};
\node[cercle] (C) at (18,11){};
\draw[-](9, -4)--(9, 2);
\draw[-](9, 2)--(3,8);
\draw[-](3,8)--(0,11);
\draw[-](5,6)--(6,8);
\draw[-](6,8)--(4, 11);
\draw[-](6,8)--(8, 11);
\draw[-](9, 2)--(13, 6);
\draw[-](13,6)--(18, 11);
\draw[-](13, 6)--(12, 8);
\draw[-](12, 8)--(10.6, 11);
\draw[-](12, 8)--(14, 11);
\path(0,11) node[above] {c};
\path(4,11) node[above] {b};
\path(8,11) node[above] {a};
\path(10.6,11) node[above] {a};
\path(14,11) node[above] {b};
\path(18,11) node[above] {d};
\end{tikzpicture}\xrightarrow[\mbox{Definition \ref{reduction}}]{\mbox{by item $8$ of }}\begin{tikzpicture}[ scale=0.1]
\tikzstyle{terminal}=[circle,draw,thick,scale=0.15mm]
\tikzstyle{cercle}=[circle,draw,thick,fill=black, scale=0.15mm]
\fill[color=black](9, -4)circle(0mm);
\fill[color=black](9, 2)circle(0mm);
\fill[color=black](3, 8)circle(0mm);
\fill[color=black](15, 8)circle(0mm);
\node[cercle] (C) at (0,11){};
\node[terminal] (T) at (6,11){};
\node[cercle] (C) at (12,11){};
\node[cercle] (C) at (18,11){};
\draw[-](9, -4)--(9, 2);
\draw[-](9, 2)--(0, 11);
\draw[-](3, 8)--(6, 11);
\draw[-](9, 2)--(18, 11);
\draw[-](15, 8)--(12, 11);
\path(0,11) node[above] {c};
\path(6,11) node[above] {\small{b+1}};
\path(12,11) node[above] {\small{b+1}};
\path(18,11) node[above] {d};
\end{tikzpicture}\xrightarrow[\mbox{Definition \ref{reduction}}]{\mbox{by item $5$ of }}\begin{tikzpicture}[ scale=0.1]
\tikzstyle{terminal}=[circle,draw,thick,scale=0.15mm]
\tikzstyle{cercle}=[circle,draw,thick,fill=black, scale=0.15mm]
\fill[color=black](8, 0)circle(0mm);
\fill[color=black](8, 6)circle(0mm);
\node[cercle] (C) at (16,15){};
\node[cercle] (C) at (0, 15){};
\draw[-](8, 0)--(8, 6);
\draw[-](8, 6)--(0, 14);
\draw[-](8, 6)--(16, 14);
\path (0,15) node[above]{c+1};
\path (16,15) node[above]{d+1};
\end{tikzpicture}
\end{align*}
\end{ex}

In general, there are different possible reductions of a Super-Leaf weighted tree $\varphi$. Nevertheless, it turns out that all possible reductions of $\varphi$ end up with the same reduced Super-Leaf weighted tree. To see this we need the following lemma.
\begin{lem}\label{reduction1}\item
For any two elementary reductions $\varphi\rightarrow\varphi_1$ and $\varphi\rightarrow\varphi_2$ of a Super-Leaf weighted tree $\varphi$, there exist elementary reductions $\varphi_1\rightarrow\varphi_0$ and $\varphi_2\rightarrow\varphi_0.$
\end{lem}
\begin{prf}\item
Let $\varphi\xrightarrow[]{\lambda_1}\varphi_1$ and  $\varphi\xrightarrow[]{\lambda_2}\varphi_2$ be two elementary reductions of a Super-Leaf weighted tree $\varphi$. There are two possible ways to carry out the reductions $\lambda_1$ and $\lambda_2.$
\begin{description}
  \item[Case1: Disjoint reductions] In this case $\varphi=\psi_1\vee\psi_2$ and $\lambda_1(\varphi)=\overline{\psi_1}\vee\psi_2$, $\lambda_2(\varphi)=\psi_1\vee\overline{\psi_2}$. Then
  \begin{align*}
    \varphi\xrightarrow[]{\lambda_1}\overline{\psi_1}\vee\psi_2\xrightarrow[]{\lambda_2}\overline{\psi_1}\vee\overline{\psi_2} \\
    \varphi\xrightarrow[]{\lambda_2}\psi_1\vee\overline{\psi_2}\xrightarrow[]{\lambda_1}\overline{\psi_1}\vee\overline{\psi_2}.
  \end{align*}
Hence the lemma holds.
  \item[Case2: Overlapping reductions] In this case $\varphi$ can take the following form
  $$\varphi=(\psi_1\vee\extredtwon)\vee(\redtwtw\vee\psi_2). $$
 Then, there are two different possible reductions of $\varphi:$
 \begin{itemize}
   \item $\varphi\xrightarrow[]{\lambda_1}\a(\psi_1)\vee(\redtwtw\vee\psi_2),$ and
   \item $\varphi\xrightarrow[]{\lambda_2}(\psi_1\vee \redtwtw)\vee\a(\psi_2).$
 \end{itemize}
To have the lemma, we would like to show that $\lambda_2(\varphi)=\lambda_1(\varphi):$
\begin{eqnarray*}
  \lambda_2(\varphi) &=&(\psi_1\vee \redtwtw)\vee\a(\psi_2)  \\
   &=&
   (\psi_1\vee \redtwtw)\vee(\mathbf{1}\vee\psi_2)  \\
  &=&
  (\psi_1\vee \redtwtw)\vee(\extredttwon\vee\psi_2)  \\
   &=&
   \a(\psi_1)\vee(\redtwtw\vee\psi_2).
\end{eqnarray*}
\end{description}
\end{prf}
\begin{pro}\label{reduction2}
Let $\varphi$ be a Super-Leaf weighted tree. Then any two reductions of $\varphi:$
\begin{align*}
  \varphi\rightarrow\varphi'_1\rightarrow\varphi'_2\rightarrow\cdots\rightarrow\varphi'_n\\
  \varphi\rightarrow\varphi"_1\rightarrow\varphi"_2\rightarrow\cdots\rightarrow\varphi"_m
\end{align*}
result in the same reduced Super-Leaf weighted tree of $\varphi,$ i. e, $\varphi'_n=\varphi"_n.$
\end{pro}
\begin{prf}\item
By induction on $\varphi.$ If $\varphi=\mathbf{1}$ or $\varphi$ is a Super-Leaf weighted $1$-tree, then $\varphi$ is reduced and there is nothing to prove.\\
Let $\varphi\in B_r,~r\geq1$ and
\begin{align*}
  \varphi\rightarrow\varphi'_1\rightarrow\varphi'_2\rightarrow\cdots\rightarrow\varphi'_n\\
  \varphi\rightarrow\varphi"_1\rightarrow\varphi"_2\rightarrow\cdots\rightarrow\varphi"_m
\end{align*}
be two reductions of $\varphi.$ Then by Lemma \ref{reduction1}, there are elementary reductions $\varphi'_1\rightarrow\varphi_0$ and $\varphi"_1\rightarrow\varphi_0.$ Consider a reduction process for $\varphi_0:$
\begin{align*}
  \varphi_0\rightarrow\varphi_1\rightarrow\varphi_2\rightarrow\cdots\rightarrow\varphi_k.
\end{align*}
By induction all reduced Super-Leaf weighted trees of $\varphi'_1$ are equal to each other, as well as all reduced Super-Leaf weighted trees of $\varphi"_1$. Since $\varphi_k$ is a reduced Super-Leaf weighted tree of both $\varphi'_1$ and $\varphi"_1$, then $\varphi'_n=\varphi_k=\varphi"_m$ as desired. This proves the proposition.
\end{prf}
Let $G$ be the set of all reduced Super-Leaf weighted trees.\\
We now intend to define inverse of a Super-Leaf weighted $n$-trees to hold the structure of  Hom-group on $G$. We first define Mirror inverse type of a Super-Leaf weighted $n$-tree($n\geq2$). \\
 \begin{df}
   Mirror inverse type of a Super-Leaf weighted $n$-trees($n\geq2$) $(\varphi,V)$ is another Super-Leaf weighted  $n$-trees $M(\varphi)$ with left and right children of all nodes interchanged and inversing the type of every leaf. We note $\mathbf{MRinv}$ the application that sends each Super-Leaf weighted $n$-trees to its mirror inverse type.
 \end{df}
 \begin{exs}
 \begin{enumerate}
   \item $$\mathbf{MRinv}\Big(\begin{tikzpicture}[ scale=0.05]
\tikzstyle{terminal}=[circle,draw,thick,scale=0.15mm]
\tikzstyle{cercle}=[circle,draw,thick,fill=black, scale=0.15mm]
\fill[color=black](8, 0)circle(0mm);
\fill[color=black](8, 6)circle(0mm);
\node[cercle] (C) at (16,15){};
\node[terminal] (T) at (1,15){};
\draw[-](8, 0)--(8, 6);
\draw[-](8, 6)--(1, 14);
\draw[-](8, 6)--(16, 14);
\path (1,15) node[above]{$7$};
\path (16,15) node[above]{$3$};
\end{tikzpicture}\Big)=\begin{tikzpicture}[ scale=0.05]
\tikzstyle{terminal}=[circle,draw,thick,scale=0.15mm]
\tikzstyle{cercle}=[circle,draw,thick,fill=black, scale=0.15mm]
\fill[color=black](8, 0)circle(0mm);
\fill[color=black](8, 6)circle(0mm);
\node[cercle] (T) at (16,15){};
\node[terminal] (C) at (1,15){};
\draw[-](8, 0)--(8, 6);
\draw[-](8, 6)--(1, 14);
\draw[-](8, 6)--(16, 14);
\path (1,15) node[above]{$3$};
\path (16,15) node[above]{$7$};
\end{tikzpicture}$$
\item $$\mathbf{MRinv}\Big(\begin{tikzpicture}[ scale=0.05]
\tikzstyle{terminal}=[circle,draw,thick,scale=0.15mm]
\tikzstyle{cercle}=[circle,draw,thick,fill=black, scale=0.15mm]
\fill[color=black](8, 0)circle(0mm);
\fill[color=black](8, 6)circle(0mm);
\fill[color=black](4, 10)circle(0mm);
\node[cercle] (C) at (0,15){};
\node[cercle] (C) at (8,15){};
\node[terminal] (T) at (16, 15){};
\draw[-](8, 0)--(8, 6);
\draw[-](8, 6)--(0, 14);
\draw[-](4, 10)--(8, 14);
\draw[-](8, 6)--(16, 14);
\path(-2,15) node [above] {$-1$};
\path(8, 15) node[above] {$1$};
\path(16, 15) node[above] {$3$};
\end{tikzpicture}\Big)=\begin{tikzpicture}[ scale=0.05]
\tikzstyle{terminal}=[circle,draw,thick,scale=0.15mm]
\tikzstyle{cercle}=[circle,draw,thick,fill=black, scale=0.15mm]
\fill[color=black](8, 0)circle(0mm);
\fill[color=black](8, 6)circle(0mm);
\fill[color=black](4, 10)circle(0mm);
\node[cercle] (C) at (0,15){};
\node[terminal] (T) at (8,15){};
\node[terminal] (T) at (16, 15){};
\draw[-](8, 0)--(8, 6);
\draw[-](8, 6)--(0, 14) ;
\draw[-](8, 6)--(17, 14) ;
\draw[-](12, 10)--(8, 14);
\path(0,15) node [above] {$3$};
\path(8, 15) node[above] {$1$};
\path(17, 15) node[above] {$-1$};
\end{tikzpicture}$$
   \item $$\mathbf{MRinv}\Big(\begin{tikzpicture}[ scale=0.1]
\tikzstyle{terminal}=[circle,draw,thick,scale=0.15mm]
\tikzstyle{cercle}=[circle,draw,thick,fill=black, scale=0.15mm]
\fill[color=black](9, -4)circle(0mm);
\fill[color=black](9, 2)circle(0mm);
\fill[color=black](4.5, 6)circle(0mm);
\fill[color=black](6, 8)circle(0mm);
\fill[color=black](15, 8)circle(0mm);
\node[terminal] (T) at (0,11){};
\node[terminal] (T) at (4,11){};
\node[cercle] (C) at (8,11){};
\node[cercle] (C) at (12,11){};
\node[cercle] (C) at (18,11){};
\draw[-](9, -4)--(9, 2);
\draw[-](9, 2)--(0, 11) ;
\draw[-](9, 2)--(18, 11);
\draw[-](4.5, 6)--(6, 8);
\draw[-](6, 8)--(4, 11);
\draw[-](6, 8)--(8, 11) ;
\draw[-](15, 8)--(12, 11);
\path(0,11) node [above] {$2$};
\path(4, 11) node[above] {$-3$};
\path(8, 11) node[above] {$6$};
\path(12, 11) node[above] {$2$};
\path(18, 11) node[above] {$4$};
\end{tikzpicture}\Big)=\begin{tikzpicture}[ scale=0.1]
\tikzstyle{terminal}=[circle,draw,thick,scale=0.15mm]
\tikzstyle{cercle}=[circle,draw,thick,fill=black, scale=0.15mm]
\fill[color=black](9, -4)circle(0mm);
\fill[color=black](9, 2)circle(0mm);
\fill[color=black](4.5, 6)circle(0mm);
\fill[color=black](6, 8)circle(0mm);
\fill[color=black](15, 8)circle(0mm);
\node[terminal] (T) at (0,11){};
\node[terminal] (T) at (6,11){};
\node[terminal] (T) at (10.6,11){};
\node[cercle] (C) at (14,11){};
\node[cercle] (C) at (18,11){};
\draw[-](9, -4)--(9, 2);
\draw[-](9, 2)--(0, 11);
\draw[-](3, 8)--(6, 11);
\draw[-](9, 2)--(18, 11);
\draw[-](9, 2)--(13, 6);
\draw[-](13, 6)--(12, 8);
\draw[-](12, 8)--(10.6, 11) ;
\draw[-](12, 8)--(14, 11);
\path(0,11) node [above] {$4$};
\path(6, 11) node[above] {$2$};
\path(10.6, 11) node[above] {$6$};
\path(14, 11) node[above] {$-3$};
\path(18, 11) node[above] {$2$};
\end{tikzpicture}$$
 \end{enumerate}
 \end{exs}
 Now, we define the inverse of  Super-Leaf weighted $n$-trees by the map $\digamma:\mathbb{T}\rightarrow \mathbb{T} $ sending $\mathbf{1}$ to $\mathbf{1}$, a Super-Leaf weighted $\Omega$-type $1$-tree $\redoneone$ to $\redonetw$, a Super-Leaf weighted $\Omega$-type $1$-tree $\redonetw$ to $\redoneone$  and a super-Leaf weighted $n$-tree ($n\geq2$) $(\varphi,a)$ to it's Mirror inverse.\\
 According to the definition \ref{reduction}, The proof of the following lemma is trivial:
 \begin{lem}
 For all $\varphi$ and $\psi$ in $\mathbb{T} $, we have
   $$\digamma(\varphi\vee\psi)=\digamma(\psi)\vee\digamma(\varphi).$$
 \end{lem}
 \begin{thm}
 Let $\widetilde{\a}$ be the bijevtive map defined on $G$ by, for all $\varphi\in\mathbb{T},~\widetilde{\a}(\overline{\varphi})=\overline{\a(\varphi)}$ and $\widetilde{\a}^{-1}(\overline{\varphi})=\overline{\a^{-1}(\varphi)}$. If we define a multiplication "$\cdot$" on $G$ by
 $$\overline{\varphi}\cdot\overline{\psi}=\overline{\varphi\vee\psi},\;\;\forall \overline{\varphi},~\overline{\psi}\in G,$$
 we have that $(G,\cdot,\mathbf{1},\widetilde{\a})$ is a regular Hom-group.\\
 $G$ will be said the free regular Hom-group with $1$-generator.
 \end{thm}
 \begin{prf}\item
 \begin{enumerate}
   \item The map $\widetilde{\a}$ is multiplicative:
   \begin{align*}
   \widetilde{\a}( \overline{\varphi}\cdot\overline{\psi})= \widetilde{\a}(\overline{\varphi\vee\psi})=\overline{\a(\varphi\vee\psi)}=\overline{\a(\varphi)\vee\a(\psi)}=\overline{\a(\varphi)}\cdot\overline{\a(\psi)}=\widetilde{\a}(\overline{\varphi})\cdot\widetilde{\a}(\overline{\psi})
   \end{align*}
   \item The multiplication "$\cdot$" is Hom-associative:
   $$\widetilde{\a}(\overline{\varphi})\cdot (\overline{\psi}\cdot\overline{\phi})=(\overline{\varphi}\cdot\overline{\psi} )\cdot \widetilde{\a}(\overline{\phi}),$$
   for any $\overline{\varphi},~\overline{\psi},~\overline{\phi}\in G$. To show this, it enough to prove that
   $$\overline{\overline{\a(\varphi)}\vee(\overline{\psi\vee\phi})}=\overline{(\overline{\varphi\vee\psi})\vee\overline{\a(\phi)}},$$
   for given $\varphi,~\psi,~\phi\in\mathbb{T}$. We will observe, that each of the reduced Super-Leaf weighted trees $\overline{\overline{\a(\varphi)}\vee(\overline{\psi\vee\phi})},$ $\overline{(\overline{\varphi\vee\psi})\vee\overline{\a(\phi)}}$ can be obtained from a same Super-Leaf weighted tree by a sequence of elementary reductions.  We will show this by induction on $\psi$.\\
   If $\psi=\mathbf{1}$, there is nothing to prove.\\
   If $\psi=\redoneone$ or $\redonetw$, we have
   \begin{equation*}
   \a(\varphi)\vee(\psi\vee\phi)=\a(\varphi)\vee(\redoneone\vee\phi)=(\varphi\vee\mathbf{1})\vee(\redoneone\vee\phi)=(\varphi\vee\begin{tikzpicture}[ scale=0.03]
\tikzstyle{terminal}=[circle,draw,thick,scale=0.15mm]
\tikzstyle{cercle}=[circle,draw,thick,fill=black, scale=0.15mm]
\fill[color=black](8, 0)circle(0mm);
\fill[color=black](8, 6)circle(0mm);
\node[cercle] (C) at (0,15){};
\node[terminal] (T) at (16, 15){};
\draw[-](8, 0)--(8, 6);
\draw[-](8, 6)--(0, 14);
\draw[-](8, 6)--(16, 14);
\path (-4,15) node[above]{$a-1$};
\path (26,15) node[above]{$a-1$};
\end{tikzpicture})(\redoneone\vee\phi)=(\varphi\vee\redoneone)\vee\a(\phi).
   \end{equation*}
If $\psi\in B_n$, then $\psi$ can take the following form $\psi=\psi_1\vee\psi_2$, where $\psi_1\in B_p$ and $\psi_2\in B_q,$ such that $p+q=n.$ By Proposition         \ref{reduction2}, we have $\overline{\psi}=\overline{\psi_1\vee\psi_2}=\overline{\psi_1}\cdot\overline{\psi_2}$. Then,
\begin{align*}
\widetilde{\a}(\overline{\varphi})\cdot(\overline{\psi}\cdot\overline{\phi})&=\widetilde{\a}(\overline{\varphi})\cdot((\overline{\psi_1}\cdot\overline{\psi_2})\cdot\overline{\phi})\\
&=
\widetilde{\a}(\overline{\varphi})\cdot(\widetilde{\a}(\overline{\psi_1})\cdot(\overline{\psi_2}\cdot\widetilde{\a}^{-1}(\overline{\phi}))\\
&=
(\overline{\varphi}\cdot \widetilde{\a}(\overline{\psi_1}))\cdot (\widetilde{\a}(\overline{\psi_2})\cdot\overline{\phi}),
\end{align*} and
\begin{align*}
 (\overline{\varphi}\cdot\overline{\psi} )\cdot \widetilde{\a}(\overline{\phi})&=(\overline{\varphi}\cdot(\overline{\psi_1}\cdot(\overline{\psi_2}) )\cdot \widetilde{\a}(\overline{\phi})\\
 &=
 ((\widetilde{\a}^{-1}(\overline{\varphi})\cdot\overline{\psi_1})\cdot \widetilde{\a}(\overline{\psi_2}))\cdot \widetilde{\a}(\overline{\phi})\\
 &=
 (\overline{\varphi}\cdot \widetilde{\a}(\overline{\psi_1}))\cdot (\widetilde{\a}(\overline{\psi_2})\cdot\overline{\phi}).
\end{align*}
Hence by Proposition \ref{reduction2}, we have  $\widetilde{\a}(\overline{\varphi})\cdot(\overline{\psi}\cdot\overline{\phi})= (\overline{\varphi}\cdot\overline{\psi} )\cdot \widetilde{\a}(\overline{\phi}).$
   \item The inverse for every $\overline{\varphi}\in G$, is defined by $\overline{\digamma(\varphi)}:$
   \begin{align*}
   \overline{\varphi}\cdot \overline{\digamma(\varphi)}=\overline{\varphi\vee\digamma(\varphi)}=\overline{\mathbf{1}}=\mathbf{1}.
   \end{align*}
 \end{enumerate}
 \end{prf}
\subsection{Free regular Hom-group}
Let $X$ be a nonempty set. We define $\mathbb{T}^X$ by $\mathbb{T}^X=\mathbf{1\oplus}\bigoplus_{n\geq1}B_n\otimes X$. Elements of $B_n\otimes X$ shall be pictured for every $\varphi\in B_n,~x_1,\cdots,x_n\in X$, by inserting, for all $i=1,\cdots,n$
the element $x_i$ at the top of the leaf with label $x_i$:
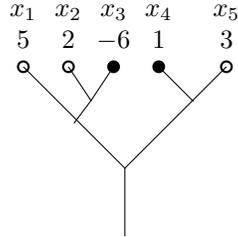
$$\begin{tikzpicture}[ scale=0.15]
\tikzstyle{terminal}=[circle,draw,thick,scale=0.15mm]
\tikzstyle{cercle}=[circle,draw,thick,fill=black, scale=0.15mm]
\fill[color=black](9, -4)circle(0mm);
\fill[color=black](9, 2)circle(0mm);
\fill[color=black](4.5, 6)circle(0mm);
\fill[color=black](6, 8)circle(0mm);
\fill[color=black](15, 8)circle(1mm);
\node[terminal] (T) at (0,11){};
\node[terminal] (T) at (4,11){};
\node[cercle] (C) at (8,11){};
\node[cercle] (C) at (12,11){};
\node[terminal] (T) at (18,11){};
\draw[-](9, -4)--(9, 2);
\draw[-](9, 2)--(0, 11);
\draw[-](9, 2)--(18, 11);
\draw[-](4.5, 6)--(6, 8);
\draw[-](6, 8)--(4, 11) ;
\draw[-](6, 8)--(8, 11) ;
\draw[-](15, 8)--(12, 11);
 \path(0,11) node [above]{$\begin{matrix}
                             x_1 \\
                             5
                           \end{matrix}$};
 \path(18,11) node[above] { $\begin{matrix}
                             x_5 \\
                             3
                           \end{matrix}$};
 \path(4,11) node [above]{ $\begin{matrix}
                             x_2 \\
                             2
                           \end{matrix}$};
 \path(8,11) node[above] { $\begin{matrix}
                            x_3 \\
                             -6
                           \end{matrix}$};
 \path(12,11) node[above] { $\begin{matrix}
                             x_4 \\
                             1
                           \end{matrix}$};
\end{tikzpicture}.$$
Let now $X^{-1}=\{x^{-1}/x\in X\}$ and consider the set $X^{*}=X\cup X^{-1}$. An element $y\in X^{*}$ can be defined as element of $B_1\otimes X$ in the following way:
\begin{enumerate}
  \item If $y=x\in X,$ we asoociate to $y$ a Super-Leaf weighted $1$-tree;
   $$y=
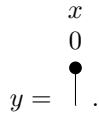
\begin{tikzpicture}[ scale=0.05]
\tikzstyle{terminal}=[circle,draw,thick,scale=0.15mm]
\tikzstyle{cercle}=[circle,draw,thick,fill=black, scale=0.15mm]
\fill[color=black](1, 0)circle(0mm);
\node[cercle] (C) at (1,10){};
\draw[-](1, 0)--(1, 10);
\path(1,10) node[above] { $\begin{matrix}
                             x \\
                             0
                           \end{matrix}$};;
\end{tikzpicture}.$$
\item If $y=x^{-1}\in X^{-1}$, for $x\in X$,    we asoociate to $y$ a Super-Leaf weighted $1$-tree;
$$y=
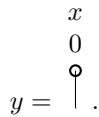
\begin{tikzpicture}[ scale=0.05]
\tikzstyle{terminal}=[circle,draw,thick,scale=0.15mm]
\tikzstyle{cercle}=[circle,draw,thick,fill=black, scale=0.15mm]
\fill[color=black](1, 0)circle(0mm);
\node[terminal] (T) at (1,10){};
\draw[-](1, 0)--(1, 10);
\path(1,10) node[above] { $\begin{matrix}
                             x \\
                             0
                           \end{matrix}$};;
\end{tikzpicture}.$$
\end{enumerate}
The operations $\vee$ and $\a$  have natural extensions to $\mathbb{T}^X$, that we denote by  same
symbols.\\
For $\a$, it defines a map $\a:X^{*}\rightarrow X^{*}$ as follows
\begin{enumerate}
  \item If $y=x\in X,$ we can define $\a(y)$ by  $$\a(y)=\begin{tikzpicture}[ scale=0.05]
\tikzstyle{terminal}=[circle,draw,thick,scale=0.15mm]
\tikzstyle{cercle}=[circle,draw,thick,fill=black, scale=0.15mm]
\fill[color=black](1, 0)circle(0mm);
\node[cercle] (C) at (1,10){};
\draw[-](1, 0)--(1, 10);
\path(1,10) node[above] { $\begin{matrix}
                             x \\
                             1
                           \end{matrix}$};
\end{tikzpicture}.$$
  \item If $y=x^{-1}\in X^{-1}$, for $x\in X$, we can define $\a(y)$ by  $$\a(y)=\begin{tikzpicture}[ scale=0.05]
\tikzstyle{terminal}=[circle,draw,thick,scale=0.15mm]
\tikzstyle{cercle}=[circle,draw,thick,fill=black, scale=0.15mm]
\fill[color=black](1, 0)circle(0mm);
\node[terminal] (T) at (1,10){};
\draw[-](1, 0)--(1, 10);
\path(1,10) node[above] { $\begin{matrix}
                             x \\
                             1
                           \end{matrix}$};
\end{tikzpicture}.$$
\end{enumerate}
Moreover, $\a^{a}(y)$, for all $y\in X^{*}$ and $a\in\Z$ is defined as follows.
\begin{enumerate}
  \item If $y=x\in X,$ we can define $\a(y)$ by  $$\a^a(y)=\begin{tikzpicture}[ scale=0.05]
\tikzstyle{terminal}=[circle,draw,thick,scale=0.15mm]
\tikzstyle{cercle}=[circle,draw,thick,fill=black, scale=0.15mm]
\fill[color=black](1, 0)circle(0mm);
\node[cercle] (C) at (1,10){};
\draw[-](1, 0)--(1, 10);
\path(1,10) node[above] { $\begin{matrix}
                             x \\
                             a
                           \end{matrix}$};
\end{tikzpicture}.$$
  \item If $y=x^{-1}\in X^{-1}$, for $x\in X$, we can define $\a(y)$ by  $$\a^a(y)=\begin{tikzpicture}[ scale=0.05]
\tikzstyle{terminal}=[circle,draw,thick,scale=0.15mm]
\tikzstyle{cercle}=[circle,draw,thick,fill=black, scale=0.15mm]
\fill[color=black](1, 0)circle(0mm);
\node[terminal] (T) at (1,10){};
\draw[-](1, 0)--(1, 10);
\path(1,10) node[above] { $\begin{matrix}
                             x \\
                             a
                           \end{matrix}$};
\end{tikzpicture}.$$
\end{enumerate}
So, the map $\a$ sends each element $\varphi\otimes (x_1,\cdots,x_n)$ of $B_n\otimes X$ to  $\a(\varphi)\otimes (x_1,\cdots,x_n).$\\
Also the reduction process and the map $\mathbf{MRinv}$ have natural extensions to $\mathbb{T}^X$. They are defined as follows.
\begin{df}
The reduction of an element of $\mathbb{T}^X$ is defined by the same process as in Definition \ref{reduction} but for two leaves labeled by the same $x\in X.$
\end{df}
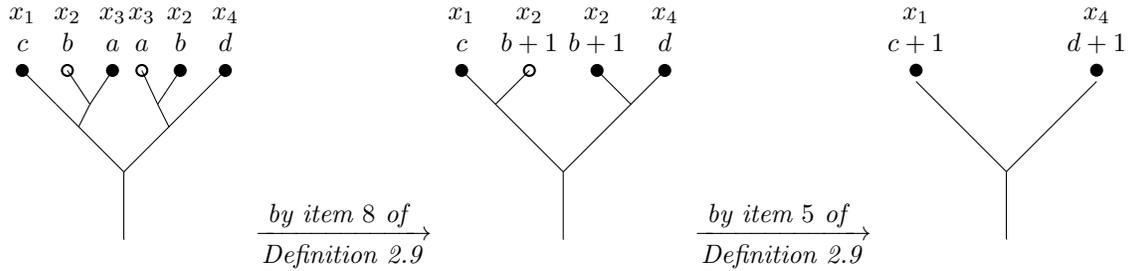
\begin{ex}
  \begin{align*}
\begin{tikzpicture}[ scale=0.15]
\tikzstyle{terminal}=[circle,draw,thick,scale=0.15mm]
\tikzstyle{cercle}=[circle,draw,thick,fill=black, scale=0.15mm]
\fill[color=black](9, -4)circle(0mm);
\fill[color=black](9, 2)circle(0mm);
\fill[color=black](6, 8)circle(0mm);
\fill[color=black](4, 11)circle(0mm);
\fill[color=black](5,6)circle(0mm);
\fill[color=black](13, 6)circle(0mm);
\fill[color=black](12, 8)circle(0mm);
\node[cercle] (C) at (0,11){};
\node[terminal] (T) at (4,11){};
\node[cercle] (C) at (8,11){};
\node[terminal] (T) at (10.6,11){};
\node[cercle] (C) at (14,11){};
\node[cercle] (C) at (18,11){};
\draw[-](9, -4)--(9, 2);
\draw[-](9, 2)--(3,8);
\draw[-](3,8)--(0,11);
\draw[-](5,6)--(6,8);
\draw[-](6,8)--(4, 11);
\draw[-](6,8)--(8, 11);
\draw[-](9, 2)--(13, 6);
\draw[-](13,6)--(18, 11);
\draw[-](13, 6)--(12, 8);
\draw[-](12, 8)--(10.6, 11);
\draw[-](12, 8)--(14, 11);
\path(0,11) node[above] { $\begin{matrix}
                             x_1 \\
                             c
                           \end{matrix}$};
\path(4,11) node[above] { $\begin{matrix}
                             x_2 \\
                             b
                           \end{matrix}$};
\path(8,11) node[above] { $\begin{matrix}
                             x_3 \\
                             a
                           \end{matrix}$};
\path(10.6,11) node[above] { $\begin{matrix}
                             x_3 \\
                             a
                           \end{matrix}$};
\path(14,11) node[above] { $\begin{matrix}
                             x_2 \\
                             b
                           \end{matrix}$};
\path(18,11) node[above] { $\begin{matrix}
                             x_4 \\
                             d
                           \end{matrix}$};
\end{tikzpicture}\xrightarrow[\mbox{Definition \ref{reduction}}]{\mbox{by item $8$ of }}\begin{tikzpicture}[ scale=0.15]
\tikzstyle{terminal}=[circle,draw,thick,scale=0.15mm]
\tikzstyle{cercle}=[circle,draw,thick,fill=black, scale=0.15mm]
\fill[color=black](9, -4)circle(0mm);
\fill[color=black](9, 2)circle(0mm);
\fill[color=black](3, 8)circle(0mm);
\fill[color=black](15, 8)circle(0mm);
\node[cercle] (C) at (0,11){};
\node[terminal] (T) at (6,11){};
\node[cercle] (C) at (12,11){};
\node[cercle] (C) at (18,11){};
\draw[-](9, -4)--(9, 2);
\draw[-](9, 2)--(0, 11);
\draw[-](3, 8)--(6, 11);
\draw[-](9, 2)--(18, 11);
\draw[-](15, 8)--(12, 11);
\path(0,11) node[above] { $\begin{matrix}
                             x_1 \\
                             c
                           \end{matrix}$};
\path(6,11)  node[above] { $\begin{matrix}
                             x_2 \\
                             b+1
                           \end{matrix}$};
\path(12,11)  node[above] { $\begin{matrix}
                             x_2 \\
                             b+1
                           \end{matrix}$};
\path(18,11) node[above] { $\begin{matrix}
                             x_4 \\
                             d
                           \end{matrix}$};
\end{tikzpicture}\xrightarrow[\mbox{Definition \ref{reduction}}]{\mbox{by item $5$ of }}\begin{tikzpicture}[ scale=0.15]
\tikzstyle{terminal}=[circle,draw,thick,scale=0.15mm]
\tikzstyle{cercle}=[circle,draw,thick,fill=black, scale=0.15mm]
\fill[color=black](8, 0)circle(0mm);
\fill[color=black](8, 6)circle(0mm);
\node[cercle] (C) at (16,15){};
\node[cercle] (C) at (0, 15){};
\draw[-](8, 0)--(8, 6);
\draw[-](8, 6)--(0, 14);
\draw[-](8, 6)--(16, 14);
\path (0,15) node[above] { $\begin{matrix}
                             x_1 \\
                             c+1
                           \end{matrix}$};
\path (16,15) node[above] { $\begin{matrix}
                             x_4 \\
                             d+1
                           \end{matrix}$};
\end{tikzpicture}
\end{align*}
\end{ex}

\begin{df}
Let $\varphi\otimes (x_1,\cdots,x_n)$ be an element of $B_n\otimes X$. Then $\mathbf{MRinv}(\varphi\otimes (x_1,\cdots,x_n))=\mathbf{Minv}(\varphi)\otimes(x_n,\cdots,x_1). $
\end{df}
\begin{ex}
 $$\mathbf{Minv}\Big(\begin{tikzpicture}[ scale=0.05]
\tikzstyle{terminal}=[circle,draw,thick,scale=0.15mm]
\tikzstyle{cercle}=[circle,draw,thick,fill=black, scale=0.15mm]
\fill[color=black](8, 0)circle(0mm);
\fill[color=black](8, 6)circle(0mm);
\fill[color=black](4, 10)circle(0mm);
\node[cercle] (C) at (0,14){};
\node[cercle] (C) at (8,14){};
\node[terminal] (T) at (16,14){};
\draw[-](8, 0)--(8, 6);
\draw[-](8, 6)--(0, 14);
\draw[-](4, 10)--(8, 14);
\draw[-](8, 6)--(16, 14);
\path(0,14) node[above] { $\begin{matrix}
                             x_1 \\
                             -1
                           \end{matrix}$};
\path(8,14) node[above] { $\begin{matrix}
                             x_2 \\
                             1
                           \end{matrix}$};
\path(16,14)node[above] { $\begin{matrix}
                             x_3 \\
                             3
                           \end{matrix}$};
\end{tikzpicture}\Big)=\begin{tikzpicture}[ scale=0.05]
\tikzstyle{terminal}=[circle,draw,thick,scale=0.15mm]
\tikzstyle{cercle}=[circle,draw,thick,fill=black, scale=0.15mm]
\fill[color=black](8, 0)circle(0mm);
\fill[color=black](8, 6)circle(0mm);
\fill[color=black](12, 10)circle(0mm);
\node[cercle] (C) at (0,14){};
\node[terminal] (T) at (8,14){};
\node[terminal] (T) at (17,14){};
\draw[-](8, 0)--(8, 6);
\draw[-](8, 6)--(0, 14) ;
\draw[-](8, 6)--(17, 14) ;
\draw[-](12, 10)--(8, 14) ;
\path(0,14) node[above] { $\begin{matrix}
                             x_3 \\
                             3
                           \end{matrix}$};
\path(8,14) node[above] { $\begin{matrix}
                             x_2 \\
                             1
                           \end{matrix}$};
\path(17,14)node[above] { $\begin{matrix}
                             x_1 \\
                             -1
                           \end{matrix}$};
\end{tikzpicture}.$$
\end{ex}
 \begin{df}
The set of all reduced elements of $B_n\otimes X$ is the free regular Hom-group.
 \end{df}
\section{Normal Hom-subgroups}
We introduce  in this section Normal Ho-subgroups and discuss their properties.
\subsection{Definitions and Examples}
\begin{df}
Let $(G,\mu_G,e,\alpha)$ be a regular Hom-group.  A Hom-subgroup $H$ of $G$ is called a normal Hom-subgroup of $G$ if
\begin{equation}\label{AB1}
  \forall g\in G, \;\;gH=Hg.
\end{equation}
We indicate that $H$ is a normal Hom-subgroup of $G$ with the notation $H\unrhd G.$
\end{df}
\begin{lem}
Let $H$ be a Hom-subgroup of a regular Hom-group $(G,\mu,e,\alpha)$. $H$ is a normal Hom-subgroup if and only if
\begin{equation}\label{AB2}
  \forall g\in G,\;\forall h\in H, \;\;(gh)\a(g^{-1})\in H.
\end{equation}
\end{lem}
\begin{prf}
   Let $H$ be a normal Hom-subgroup of a regular Hom-group $(G,\mu_G,e,\alpha).$ Then, for all $g\in G$ and $h\in H$, we have $gh=hg$. Multiply this equation on the right by $\a(g^{-1})$, we obtain
        \begin{align*}
        (gh)\a(g^{-1})=(hg)\a(g^{-1})=\a(h)(g g^{-1})=\a^2(h)\in H.
        \end{align*}
        So,  for all $g\in G$ and $h\in H$, we have  $(gh)\a(g^{-1})=\a(g)(hg^{-1})\in H$. Conversely, let $g\in G$ and $h\in H$ such that $(gh)\a(g^{-1})\in H$. There exists $h'\in H$ such that $(gh)\a(g^{-1})=\a^2(h')$. By
multiplying this last equality on the right by $\a^2(g)$, we find
\begin{equation*}
\Big((gh)\a(g^{-1})\Big)\a^2(g)=\a^2(h')\a^2(g),
\end{equation*}
which gives us
\begin{equation*}
\a(gh)\a(g^{-1}g)=\a^2(h'g),
\end{equation*}
or
\begin{equation*}
\a^2(gh)=\a^2(h'g).
\end{equation*}
So, the invrtibility of the map $\a$ gives that $gh=h'g$. Therefore, we have $gH=Hg$.
\end{prf}
The following lemmas give examples of normal subgroups.
\begin{lem}
Let $(G,\mu,e,\alpha)$ be a regular Hom-group and $Z(G)=\{g\in G : gx=xg\;\forall x\in G\}$ be the center of the Hom-group $G$. Then $Z(G)$ is a normal Hom-subgroup.
\end{lem}
\begin{prf}
Let $g$ be an element in $G$ and $h\in Z(G).$ Then $(gh)\a(g^{-1})=(hg)\a(g^{-1})=\a(h)(gg^{-1})=\a(h)e=\a^2(h).$ According to Lemma \ref{Z}, one can show that $(gh)\a(g^{-1})\in Z(G).$
\end{prf}
\begin{df}
Let $(G,\cdot,e,\alpha)$ be a regular Hom-group and $H\leq G.$ The centralizer $C_G(H)$ and the normalizer $N_G(H)$ of the Hom-subgroup $H$ are defined respectively as follows
\begin{itemize}
  \item $C_G(H)=\{g\in G : \forall h\in H,~(gh)\a(g^{-1})=\a^2(h)\},$
  \item $N_G(H)=\{g\in G : \forall h\in H,~(gh)\a(g^{-1})\in H\}.$
\end{itemize}
\end{df}
\begin{lem}
$C_G(H)$ and $N_G(H)$ are two Hom-subgroups of $G$ and $C_G(H)\unrhd N_G(H).$
\end{lem}
\begin{prf}
Let us show  first that $C_G(H)\leq G;$
\begin{itemize}
  \item For all $h\in H$, we have $(e\cdot h)\a(e)=\a(h)e=\a^2(h)$. Then $e\in C_G(H).$
  \item Let $g,g'\in C_G(H).$ We want to prove that $gg'\in C_G(H);$
  \begin{eqnarray*}
    ((gg')h)\a(g'^{-1}g^{-1}) &=& (\a(g)(g'\a^{-1}(h)))\a(g'^{-1}g^{-1}) 
     =
     \a^2(g)((g'\a^{-1}(h))(g'^{-1}g^{-1}))  \\
     &=&
     \a^2(g)((( \a^{-1}(g')\a^{-2}(h))g'^{-1})\a(g^{-1})) 
     =
     \a^2(g)(\a^{-1}((g'\a^{-1}(h))\a(g'^{-1}))\a(g^{-1}))  \\
     &=&
     \a^2(g)(h \a(g^{-1})) 
    =
     (\a(g)h) \a^2(g^{-1})
     =
     \a((g\a^{-1}(h))\a(g^{-1}))\\
     &=&
     \a^2(h).
  \end{eqnarray*}
Consequently, we have $gg'\in C_G(H).$
  \item Let $g\in C_G(H).$ We obtain
  \begin{eqnarray*}
    (g^{-1}h)\a(g) &=& (g^{-1}(\a^{-2}(gh)\a^{-1}(g^{-1})))\a(g) 
     =
     (g^{-1}(\a^{-1}(g)\a^{-2}(hg^{-1})))\a(g)\\
     &=&
     (\a^{-1}(g^{-1} g)\a^{-1}(hg^{-1}))\a(g)  
     =
      (hg^{-1})\a(g)
     =
     \a(h)(g^{-1}g)\\
     &=&
     \a^2(h).
  \end{eqnarray*}
 So, for all $g\in C_G(H),$ we have $g^{-1}\in C_G(H)$. Therefore, $C_G(H)\leq G$. Similarly, we can show that $N_G(H)\leq G.$\\
 \end{itemize}
 Now, we prove that $C_G(H)\unrhd N_G(H).$ Let $n\in N_G(H)$ and $c\in C_G(H),$ for all $h\in H$, we have
 \begin{eqnarray*}
   [((nc)\a(n^{-1}))h]\a[\a(n)(c^{-1}n^{-1}]& = &[\a(nc)(\a(n^{-1})\a^{-1}(h))]\a[\a(n)(c^{-1}n^{-1})] \\
  & =&
   \a^2(nc)[(\a(n^{-1})\a^{-1}(h))(\a(n)(c^{-1}n^{-1}))]\\
  & =&
   \a^2(nc)[((n^{-1}\a^{-2}(h))\a(n))\a(c^{-1}n^{-1})] \\  
   &=&
   \a^2(nc)(h'\a(c^{-1}n^{-1})\;\;\mbox{where}\;\;h'=(n^{-1}\a^{-2}(h))\a(n)  \\
   &=&
   \a^2(nc)((\a^{-1}(h')\a(c^{-1}))\a^2(n^{-1}))\\
   &=&
   (\a(nc)(\a^{-1}(h')\a(c^{-1})))\a^3(n^{-1})\\
   &=&
   (\a^2(n)(\a(c)(\a^{-2}(h')c^{-1})))\a^3(n^{-1})\\
   &=&
   (\a^2(n)h')\a^3(n^{-1}).
 \end{eqnarray*}
 So $(nc)\a(n^{-1})$ is in $C_G(H)$. Indeed
 \begin{eqnarray*}
   (\a^2(n)h')\a^3(n^{-1}) &=& (\a^2(n)((n^{-1}\a^{-2}(h))\a(n)))\a^3(n^{-1}) \\
    &=&
    ((\a(n)(n^{-1}\a^{-2}(h)))\a^2(n))\a^3(n^{-1}) \\
    &=&
    (h\a^2(n))\a^3(n^{-1})\\
    &=&
    \a^2(h).
 \end{eqnarray*}
\end{prf}
\begin{lem}
Every Hom-subgroup of an abelian Hom-group is normal.
\end{lem}
\begin{lem}
Let $(G\times H,\cdot,e,\a)$ be the Hom-group defined in (\ref{produit}). Then $(G\times\{e_H\},\cdot,e,\a)$  is a normal Hom-subgroup of $(G\times H,\cdot,e,\a)$
\end{lem}
\begin{prf}
For all $g_1,~g_2\in G$ and $h\in H$, we have
\begin{align*}
  ((g_1,h)(g_2,e_H))\a(g_1^{-1},h^{-1})=((g_1g_2)\a_G(g_1^{-1}),\a_H(hh^{-1}))=   ((g_1g_2)\a_G(g_1^{-1}),e_H).
\end{align*}
So, we have the lemma.
\end{prf}
\begin{pro}
Let $f:G\rightarrow H$ be a Hom-group homomorphism between the Hom-group $(G,\cdot,e_G,\a_G)$  and the Hom-group $(H,\cdot,e_H,\a_H)$.
If $f(e_G)=e_H$, then $Ker(f)=\{g\in G/f(g)=e_H\}$ is a normal Hom-subgroup of $G.$
\end{pro}
\begin{prf}
According to \cite{H1}, $Ker(f)$ is a Hom-subgroup of $G.$ Let now $g\in G$ and $h\in Ker(f)$. We have
\begin{align*}
  f((gh)\a_G(g^{-1}))=(f(g)f(h))f(\a_G(g^{-1}))=   (f(g)e_H)\a_H(f(g^{-1}))=\a_H(f(g)f(g^{-1}))=\a_H(e_H)=e_H.
\end{align*}
Therefore $Ker(f)\unrhd G.$
\end{prf}
A consequence of this proposition is the following corollary
\begin{cor}\label{simplealpha}
Let $(G,\cdot,e,\a)$ be a Hom-group. Then $Ker(\a)\unrhd G.$
\end{cor}
\begin{pro}\label{salpha}
Let $f:G\rightarrow H$ be a Hom-group homomorphism between the Hom-group $(G,\cdot,e_G,\a_G)$  and the Hom-group $(H,\cdot,e_H,\a_H)$ such that $f(e_G)=e_H$.
\begin{enumerate}
  \item If $N\unrhd G$, then $f(N)\unrhd f(G).$
  \item If $M\unrhd f(G)$, then $f^{-1}(M)\unrhd G.$
\end{enumerate}
\end{pro}
\begin{prf}
\begin{enumerate}
  \item Let $h\in f(G)$ and $n\in N$. If $g\in G$ such that $f(g)=h,$ we have $$(hf(n))\a_H(h^{-1})=(f(g)f(n))\a_H(f(g^{-1}))=f(gn)f(\a_G(g^{-1}))=f((gn)\a_G(g^{-1}))\in f(N),$$
since $N\unrhd G$. So $f(N)\unrhd f(G).$
  \item Let $g\in G$ and $n\in f^{-1}(M)$, we have
  $$f((gn)\a_G(g^{-1}))=(f(g)f(n))f(\a_G(g^{-1}))=(f(g)f(n))\a_H(f(g^{-1}))\in M, $$
since $M\unrhd f(G)$, proving the proposition.
\end{enumerate}
\end{prf}
\subsection{Quotient Hom-group}
\begin{df}
  If $H$ is a Hom-subgroup of a Hom-group $(G,\cdot,e,\alpha)$ and $g\in G$ then
  $$gH=\{gh/h\in H\},$$
  is a left coset of $H$.
\end{df}
\begin{lem}
Let $H$ be a Hom-subgroup of a regular Hom-group $(G,\cdot,e,\alpha)$ and $g,g'\in G$. Suppose $\a(H)=H.$ Then $gH=g'H$ if and only if $g^{-1}g'\in H.$
\end{lem}
\begin{prf}
Suppose $gH=g'H$. Then $\a(g')\in g'H$ and so $\a(g')\in gH$. Thus $\a(g')=gh$ for some $h\in H$ and we see that $g^{-1}g'=\a^{-1}(h)\in H$, since $\a(H)=H.$ Conversely, suppose $g^{-1}g'\in H.$ Then $g^{-1}g'=h$, for some $h\in H$. Thus $\a^2(g')=\a(g)h $ and consequently $\a^2(g')H=(\a(g)h)\a(H)=\a^2(g)(hH)$. Observe that $hH=H$ because $h\in H$.
Therefore $gH=g'H$, since $\a^2(g) H=\a^2(gh)$.
\end{prf}
\begin{thm}
Let $H$ be a Hom-subgroup of a regular Hom-group $(G,\cdot,e,\alpha)$ and consider the set $G/H=\{gH/g\in G\}$.
If $H\unrhd G$, then  $(G/H,\cdot,eH,\widetilde{\alpha})$ is a Hom-group, where
\begin{itemize}
  \item $\widetilde{\alpha}:G/H\rightarrow G/H$ is defined by $\widetilde{\alpha}(gH)=\a(g)H$,
  \item $g_1H\cdot g_2H=(g_1g_2)H.$
\end{itemize}
\end{thm}
\begin{prf}
\begin{enumerate}
  \item For all $g_1,~g_2\in G$, we have
  \begin{align*}
    (g_1H)(g_2H)&=\a(g_1)(H\a^{-1}(g_2H))=\a(g_1)(\a^{-1}(Hg_2)H)=   \a(g_1)(\a^{-1}(g_2H)H)\\
    &=
    \a(g_1)(g_2\a^{-1}(HH))
    =
    (g_1g_2)(H H)=
    (g_1g_2)H.
  \end{align*}
  \item Since the multiplication of $G$ is Hom-associative, then the multiplication of $G/H$ is Hom-associative.
  \item $eH$ is the identity in $G/H$.
  \item Then inverse of $gH$ is $g^{-1}H.$
\end{enumerate}
\end{prf}
\subsection{Commutator Hom-subgroup}
\begin{df}
Let $(G,\mu,e,\alpha)$ be a regular Hom-group. The Commutator Hom-subgroup $N=[G,G]$ of $G$ is the set of the elements
\begin{equation}\label{AB3}
  [g,h]=(g^{-1}h^{-1})(gh),\;\;\forall g,~h\in G.
\end{equation}
\end{df}
\begin{pro}
 The Commutator subgroup $N=[G,G]$ of a regular Hom-group $(G,\mu,e,\alpha)$ is a normal Hom-subgroup.
\end{pro}
\begin{prf}
Let $s,~g,~h\in G$. Then we have
\begin{align*}
 & [(sg)\a(s^{-1}),(sh)\a(s^{-1}) ]= \Big(\Big(\a(s)(g^{-1}s^{-1})\Big)\Big(\a(s)(h^{-1}s^{-1})\Big)\Big)\Big(\Big((sg)\a(s^{-1})\Big)\Big((sh)\a(s^{-1})\Big)\Big)\\
   &=
   \Big(\a^2(s)\Big((g^{-1}s^{-1})\Big(s\a^{-1}(h^{-1}s^{-1})\Big)\Big)\Big)\Big(\Big((sg)\a(s^{-1})\Big)\Big(\a(s)(hs^{-1})\Big)\Big)  \\
   &=
   \Big(\a^2(s)\Big((g^{-1}s^{-1}) \Big(\a^{-1}(sh^{-1})s^{-1}\Big)\Big)\Big(\a(sg)\Big(\a(s^{-1})(s\a^{-1}(hs^{-1}))\Big)\Big) \\
   &=
   \Big(\a^2(s)\Big(\a^{-1}((g^{-1}s^{-1})(sh^{-1}))\a(s^{-1})\Big)\Big) \Big(\a(sg)\Big((s^{-1}s)(hs^{-1})\Big)\Big) \\
   &=
   \Big(\a^2(s)\Big((g^{-1}h^{-1})\a(s^{-1})\Big)\Big)\Big(\a(sh)\a(hs^{-1})\Big)
   =
   \a^3(s)\Big(\Big((g^{-1}h^{-1})\a(s^{-1})\Big)\Big((sg)(hs^{-1})\Big)\Big)  \\
   &=
   \a^3(s)\Big(\a(g^{-1}h^{-1})\Big(\a(s^{-1})\a^{-1}\Big((sg)(hs^{-1})\Big)\Big)\Big) 
   =
   \a^3(s)\Big(\a(g^{-1}h^{-1})\Big(\Big(s^{-1}\a^{-1}(sg)\Big)(hs^{-1})\Big)\Big)  \\
  & =
   \a^3(s)\Big(\a(g^{-1}h^{-1})\Big(\Big(\a^{-1}(s^{-1}s)g\Big)(hs^{-1})\Big)\Big)   
   =
   \a^3(s)\Big(\a(g^{-1}h^{-1})\Big(\a(g)(hs^{-1})\Big) \Big)\\
   &=
   \a^3(s)\Big(\a(g^{-1}h^{-1})\Big((gh)\a(s^{-1})\Big)\Big) 
   =
   \a^3(s)\Big(\Big((g^{-1}h^{-1})(gh)\Big)\a^2(s^{-1})\Big)  
   =
   \a^3(s)\Big([g,h]\a^2(s^{-1})\Big)\\
   &=
   \Big(\a^2(s)[g,h]\Big) \a^3(s^{-1}).
\end{align*}
So, we can notice that $(s[g,h])\a(s^{-1})=[(\a^{-2}(s)g)\a^{-1}(s^{-1}),(\a^{-2}(s)h)\a^{-1}(s^{-1})]$. This completes the proof.
\end{prf}

\begin{df}
Let $(G,\mu,e,\alpha)$ be a regular Hom-group and $N=[G,G]$ be the Commutator Hom-subgroup  of $G$. We define the quotient Hom-group $(G^{ab}=G/N,\mu',\widetilde{\alpha})$ by
\begin{enumerate}
  \item $\mu'(gN,hN)=\mu(g,h) N,$
  \item $\widetilde{\a}(gN)=\a(g) N.$
\end{enumerate}
The Hom-group $G^{ab}$ is called the Abelianization of the Hom-group $G.$
\end{df}
\begin{lem}
  Let $(G,\mu,e,\alpha)$ be a regular Hom-group and $N=[G,G]$ be the Commutator subgroup  of $G$. Then the Hom-group $(G^{ab},\mu',\widetilde{\alpha})$ is abelian. Moreover, if $H$ is another normal Hom-subgroup and the quotient $G/H$ is abelian, then $N\subset H.$
\end{lem}
\begin{prf}
From Lemma \ref{RRR1}, one can show that $(G^{ab},\mu',\widetilde{\alpha})$ is abelian. The proof of the second part is straightforward.

\end{prf}
\begin{lem}
Let $(G,\mu,e,\alpha)$ be a regular Hom-group and  $(G^{ab},\mu',\widetilde{\alpha})$  its Abelianization. Then the map $\pi: G\rightarrow G^{ab}$ is an homomorphism of Hom-groups. It is called the canonical Hom-projection of $G$ onto $G^{ab}.$
\end{lem}
\begin{prf}
The proof is straightforward.
\end{prf}
\begin{pro}
Let $(G,\mu_G,e_G,\alpha_G)$ be a regular Hom-group, $(G^{ab},\mu',\widetilde{\alpha_G})$ its Abelization and $(H,\mu_H,e_H,\alpha_H)$ be a regular abelain Hom-group. If $f: G\rightarrow H$ is a Hom-group Homomorphism, then there exists a unique Hom-group homomorphism $\widetilde{f}:G^{ab}\rightarrow H$ such that $\widetilde{f}\circ\pi=f.$
\end{pro}
\begin{prf}
 We can deduce that $Ker(f)$ is a normal Hom-subgroup of $G$. Moreover, one can check that $N\subset Ker(f).$ This allows us to define a map  $\widetilde{f}:G^{ab}\rightarrow H$ by setting $\widetilde{f}(gN)=f(g).$ To see that this is well-defined, suppose that $gN=g'N$, so that $gg'^{-1}\in N\subset Ker(f)$. Then we have $f(g)f(g')^{-1}=f(gg'^{-1})=e_H$, which implies that $f(g)=f(g').$ Also note that  $\widetilde{f}:G^{ab}\rightarrow H$ is a homomorphism because
\begin{itemize}
    \item $\widetilde{f}((gg')N)=f(gg')=f(g)f(g')=\widetilde{f}(gN)\widetilde{f}(g'N),$
    \item $\widetilde{f}(\widetilde{\a_G}(gN))=\widetilde{f}(\a_G(g)N)=f(\a_G(g))=\a_H(f(g))=\a_H(\widetilde{f}(gN)),$
  \end{itemize}
and note that $\widetilde{f}\circ\pi=f$, since for all $g\in G$ we have
$$\widetilde{f}\circ\pi(g)=\widetilde{f}(gN)=f(g).$$
Finally, suppose that $F:G^{ab}\rightarrow H$ is  another morphism satisfying $F\circ\pi=f.$ Then for all $g\in G$ we have
$$\widetilde{f}(gN)=f(g)=F(\pi(g))=F(gN).$$
So that $F=\widetilde{f}$ as desired.
\end{prf}
\section{Simple Hom-groups}
\begin{df}
A Hom-group is simple if it is nontrivial and has no proper nontrivial
normal Hom-subgroups.
\end{df}
The following theorem gives a characterization of  simple Hom-groups.
\begin{thm}
Nontrivial simple Hom-groups are regular.
\end{thm}
\begin{prf}
Let $(G,\mu,e,\alpha)$ be a simple Hom-group. According to Corollary \ref{simplealpha}, we have $ker(\a)\unrhd G$. Since $(G,\mu,e,\alpha)$ is a simple Hom-group, then $Ker(\a)=\{e\}$ or $Ker(\a)= G$. The Hom-group $G$ is nontrivial, so $Ker(\a)\neq G$. Thus $\a$ is bijective.
\end{prf}
\begin{df}
Given a Hom-group $(G,\mu,e,\alpha)$, a Hom-subgroup $H$ is a maximal subgroup if:
\begin{enumerate}
  \item $H\unrhd G,$
  \item If $K\leq G$ is a normal Hom-subgroup and $H\leq K$, then $K=H$ or $K=G$, i.e., the only normal Hom-subgroup of $G$ which contains $H$ as a proper Hom-subgroup is $G.$
\end{enumerate}
\end{df}
\begin{pro}
A normal Hom-subgroup $H$ is a maximal Hom-subgroup if and only if $G/H$ is a simple Hom-group.
\end{pro}
\begin{prf}
We will start proving that if $H$ is not maximal normal, then $G/H$ is not simple. So, assume $H$ is a normal Hom-subgroup which is not maximal, that is, there is $K\unrhd G$ such that $H\lneq K\lneq G.$  Let $\pi : G\rightarrow G/H$ be the quotient map. Since $\pi$ is a Hom-group homomorphism, then $\pi(K)$ is normal. Since $H\lneq K\lneq G$, this is a nontrivial subgroup of $G/H$, showing that $G/H$ is not simple.\\
Now, we prove the converse. Assume that $G/H$ is not simple and let $\widetilde{K}$ be a nontrivial normal Hom-subgroup of $G/H$. Now consider the set $K=\pi^{-1}(\widetilde{K})$. According to Proposition \ref{salpha}, we have $K\unrhd G.$ So it remains to claim that $K$ contains $H$. Let $h\in H$, then $\a(h)\in H$. So $hH=eH$ and therefore, $ \pi(h)=eH$. Consequently, we have $h\in\pi^{-1}(eH)$ and $H\subset \pi^{-1}(eH)\subset K.$
\end{prf}
\section{Tensor product of regular Hom-groups}
\begin{df}
Let $(A,\alpha_A)$ and $(B,\alpha_B)$  be  two hom-groups and $(C,\alpha_C)$ is another Hom-group (say in additive notation even if it is not assumed to be commutative). A map $f:A\times B\rightarrow C$ from the cartesian product of any two Hom-groups $A$ and $B$ into a third $C$ is called a Hom-bilinear map if;
\begin{enumerate}
  \item $f(a_1a_2,\a_B(b))=f(a_1,b)+f(a_2,b)$,
  \item $f(\a_A(a),b_1b_2)=f(a,b_1)+f(a,b_2)$,
  \item $f(\a_A(a),\a_B(b))=\a_C\circ f(a,b)$,
\end{enumerate}
for all $a,~a_1,~a_2\in A$ and $b,~b_1,~b_2\in B$. The set of all such Hom-bilinear map from $A\times B$ to $C$ is denoted by $HL(A,B;C).$
\end{df}
\begin{rem}
  If $\a_A=\a_B=\a_C=Id,$ we have the definition of a bilinear map.
\end{rem}
\begin{lem}\label{TT1}
 Let $f$ be a Hom-bilinear map from the Cartesian product of two regular Hom-groups $(A,\alpha_A)$ and $(B,\alpha_B)$ into a third regular  Hom-group $(C,\alpha_C)$ (say in additive notation). Then $$f(e_A,b)=f(a,e_b)=0,\;\forall a\in A\;\mbox{and}\;b\in B,$$
 where $e_A$ is the neutral element of $A$ and  $e_B$ is the neutral element of $B.$
\end{lem}
\begin{prf}
For all $b\in B$, we have
$$\a_C(f(e_A,b))=f(\a_A(e_A),\a_B(b))=f(e_Ae_A,\a_B(b))=f(e_A,b)+f(e_A,b).$$
Let $c$ be the additive inverse of $f(e_A,b)$ in $C$.Then
\begin{align*}
  0 &=\a_C(c+f(e_A,b))=
   \a_C(c)+\a_C(f(e_A,b))=
   \a_C(c)+(f(e_A,b)+f(e_A,b))\\
   &=
   (c+f(e_A,b))+\a_C(f(e_A,b))=0+\a_C(f(e_A,b))=
   \a^2_C(f(e_A,b)).
\end{align*}
Since $\a_C$ is injective, we have $f(e_A,b)=0.$ The same holds for $f(a,e_b)=0.$
\end{prf}
\begin{cor}
Let $f:A\times B\rightarrow C$ be a Hom-bilinear map from the Cartesian product of two regular Hom-groups $(A,\alpha_A)$ and $(B,\alpha_B)$ into a third regular  Hom-group $(C,\alpha_C)$ (say in additive notation). Then for all $a\in A$ and $b\in B$, we have $f(a^{-1},b)=-f(a,b)$ and $f(a,b^{-1})=-f(a,b).$
\end{cor}
\begin{lem}
  Let $f$ be a Hom-bilinear map from the Cartesian product of two regular  Hom-groups $(A,\alpha_A)$ and $(B,\alpha_B)$ into a third regular Hom-group $(C,\alpha_C)$ (the product of $C$ is denoted  additively  even if it is not commutative). Then, for all  $a_1,~a_2\in A$ and $b_1,~b_2\in B$, we have
  \begin{equation}\label{T1}
    f(a_1,b_2)+f(a_2,b_1)= f(a_2,b_1)+f(a_1,b_2),
  \end{equation}
i.e., any two elements of the image of $f$ commute.
\end{lem}
\begin{prf}
  Let $a_1,~a_2\in A$ and $b_1,~b_2\in B$. Then we have
  \begin{align*}
    f(\a_A^2(a_1+a_2),\a_B^2(b_1+b_2)) &= f(\a_A^2(a_1+a_2),\a_B^2(b_1)+\a_B^2(b_2))\\
     &=
     f(\a_A(a_1+a_2),\a_B^2(b_1))+f(\a_A(a_1+a_2),\a_B^2(b_2))  \\
     &=
     \Big(f(\a_A(a_1),\a_B(b_1))+f(\a_A(a_2),\a_B(b_1))\Big)\\
     &\ \ \ +\Big(f(\a_A(a_1),\a_B(b_2))+f(\a_A(a_2),\a_B(b_2))\Big)  \\
     &=
     \a_C\Big( \Big(f(a_1,b_1)+f(a_2,b_1)\Big)+\Big(f(a_1,b_2)+f(a_2,b_2)\Big)\Big).
  \end{align*}
  On other hand we have
   \begin{align*}
    f(\a_A^2(a_1+a_2),\a_B^2(b_1+b_2)) &= f(\a_A^2(a_1)+\a_A^2(a_2)),\a_B^2(b_1+b_2))\\
     &=
     f(\a_A^2(a_1),\a_B(b_1+b_2))+f(\a_A^2(a_2),\a_B(b_1+b_2))  \\
     &=
     \Big(f(\a_A(a_1),\a_B(b_1))+f(\a_A(a_1),\a_B(b_2))\Big)\\
     &\ \ \ +\Big(f(\a_A(a_2),\a_B(b_1))+f(\a_A(a_2),\a_B(b_2))\Big)  \\
     &=
     \a_C\Big( \Big(f(a_1,b_1)+f(a_1,b_2)\Big)+\Big(f(a_2,b_1)+f(a_2,b_2)\Big)\Big).
  \end{align*}
So, it follows from Lemma \ref{RR1} that $  f(a_1,b_2)+f(a_2,b_1)= f(a_2,b_1)+f(a_1,b_2).$
\end{prf}
\begin{lem}\label{TT2}
 Let $f$ be a Hom-bilinear map from the Cartesian product of two regular  Hom-groups $(A,\alpha_A)$ and $(B,\alpha_B)$ into a third regular Hom-group $(C,\alpha_C)$ (the product of $C$ is denoted additively even if it is not commutative). Then,
 \begin{enumerate}
   \item for all $b\in B$, the kernel of the map $f(\cdot,b):a\in A\rightarrow f(a,b)\in C$ contains $[A,A]$,
   \item for all $a\in A$, the kernel of the map $f(a,\cdot):b\in B\rightarrow f(a,b)\in C$ contains $[B,B]$.
 \end{enumerate}
\end{lem}
\begin{prf}
By straightforward calculation, we have
\begin{equation*}
f([a_1,a_2],b)=-(f(a_1,\a_B^{-2}(b))+f(a_2,\a_B^{-2}(b)))+(f(a_1,\a_B^{-2}(b))+f(a_2,\a_B^{-2}(b))),\;\forall [a_1,a_2]\in [A,A], \ \ b\in B,
\end{equation*}
and
\begin{equation*}
f(a,[b_1,b_2])=-(f(\a_A^{-2}(a),b_1)+f(\a_A^{-2}(a),b_2))+(f(\a_A^{-2}(a),b_1)+f(\a_A^{-2}(a),b_2)),\;\forall [b_1,b_2]\in [B,B],\ \ a\in A.
\end{equation*}
\end{prf}
\begin{df}
 Let $(A,\alpha_A)$ and $(B,\alpha_B)$ be two regular  Hom-groups.
A tensor product  of $A$ and $B$ is a regular  Hom-group $(A\otimes B, \a_{A\otimes B})$ which is equipped with a Hom-bilinear map
  $$\t:A\times B\rightarrow A\otimes B,$$
 such that, for every regular  Hom-group $(C,\alpha_C)$ and every Hom-bilinear map $f:A\times B \rightarrow C$, there exists a unique Hom-group homomorphism $\widetilde{f}:A\otimes B \rightarrow C$ such that the diagram:\\
 \vskip 0.3cm
$$\xymatrix{ A\times B \ar[rr]^\t \ar^{f}[rd] && A\otimes B \ar[ld]^{\widetilde{f}} \\ & C },$$
commutes, that is $\widetilde{f}\circ \t=f.$
 \end{df}
 \begin{pro}
  Let $(A,\alpha_A)$ and $(B,\alpha_B)$ be two regular Hom-groups. Tensor products of $A$ and $B$ are unique up to  unique isomorphism.
 \end{pro}
 \begin{lem}
  Let $(A,\alpha_A)$ and $(B,\alpha_B)$ be two  Hom-groups.  We have $(A\otimes B, \t, \a_{A\otimes B})\cong(B\otimes A, \t', \a_{B\otimes A}).$
 \end{lem}
 \begin{prf}
 Defining the maps
 \begin{equation*}
 \begin{array}{ccccc}
f & : & A\times B& \to & B\times A\\
 & & (a,b) & \mapsto & f(a,b)=(b,a) \\
\end{array},\ \ \
\begin{array}{ccccc}
g & : & B\times A& \to & A\times B\\
 & & (b,a) & \mapsto & f(b,a)=(a,b)
\end{array},
 \end{equation*}
we have $f\circ(\a_A\times \a_B)=(\a_B\times\a_A) \circ f$ and $g\circ(\a_B\times \a_A)=(\a_A\times\a_B) \circ g.$ From  the universel property of $\t$ and $\t'$, we have the following diagram:
\begin{displaymath}
\xymatrix{
A\times B \ar[d]^\t
\ar[r]^f
& {B\times A} \ar[r]^g \ar[d]^{\t'}
& {A\times B} \ar[d]^\t\\
{A\otimes B}
\ar[r]^{f'}
& {B\otimes A} \ar[r]^{g'}
& {A\otimes B}
}
\end{displaymath}
The Hom-bilinearity  of $\t'\circ f$ and $\t\circ g$ ensure the existence of  homomorphisms $f'$ and $g'$. Now $g\circ f$ is the identity map of $A\times B$, hence $g'\circ f'$ is the identity map of $A\otimes B$. Similarly $f'\circ g'$ is the identity of $B\otimes A$, thus  $(A\otimes B, \t, \a_{A\otimes B})\cong(B\otimes A, \t', \a_{B\otimes A}).$
\end{prf}
\begin{pro}
Let $(A,\alpha_A)$ and $(B,\alpha_B)$ be two regular Hom-groups. Then $A\otimes B\cong A^{ab}\times B^{ab}$. In particular, $A\otimes B$ is a regular abelian Hom-group.
\end{pro}
\begin{prf}
According to Lemma \ref{TT2}, one can show that any Hom-bilinear map $f:A\times B\rightarrow C$ from the Cartesian product of two regular Hom-groups $(A,\alpha_A)$ and $(B,\alpha_B)$ into a third regular  Hom-group $(C,\alpha_C)$ is determined by a unique  Hom-bilinear $f':A^{ab}\times B^{ab}\rightarrow C$ such that $f=f'\circ\gamma$, where $\gamma:A\times B\rightarrow A^{ab}\times B^{ab}$ is the canonical projection map.
\end{prf}
\section{Hom-type rings }
In this section, we introduce Hom-rings and their properties. 
\subsection{Definitions and properties}
The word Hom-ring was mentioned for the first time by A. Gohr and Y. Fregier in \cite{FG}. 
They defined  Hom-rings  as follows.
Let $A$ be a set together with two binary operations $+:A\times A\rightarrow A$ and  $\cdot:A\times A\rightarrow A$, one self-map $\a:A\rightarrow A$ and a special element $0\in A$ . Then $(A,+,\cdot,0,\a)$ is called a Hom-ring if:
\begin{itemize}
  \item $(A,+,0)$ is an abelian group.
  \item The multiplication is distributive on both sides.
  \item $\a$ is an abelian group homomorphism.
  \item $\a$ and $\cdot$ satisfy the Hom-associativity condition; $\a(x)\cdot(y\cdot z)=(x\cdot y)\cdot \a(z)$.
\end{itemize}

In this section, we define different types of Hom-rings based on Hom-groups and provide various ways of constructing them.
\begin{df}\label{homring1}
A Hom-ring of type $(1)$ is a  tuple $(A,+,\cdot,0,\a,\b)$ consisting of a  set $A$ together with two binary operations $+:A\times A\rightarrow A$  (the addition) and  $\cdot:A\times A\rightarrow A$ (the multiplication) and two set maps  $\a,\b :A\rightarrow A$, such that:
\begin{enumerate}
  \item $(A,+,0,\a)$ is an abelian Hom-group.
  \item $\b$ is an endomorphism of the abelian Hom-group $(A,+,0,\a)$, i.e, $\b(x+y)=\b(x)+\b(y),\;\forall x,~y\in A,$ and $\a\circ\b=\b\circ\a.$
  \item $\a$ and $\b$ are multiplicative maps , i.e, $\a(xy)=\a(x)\a(y)$ and  $\b(xy)=\b(x)\b(y)$, for all $x,~y\in A.$
  \item The map $\b$ and the product $\cdot$ satisfy the Hom-associativity condition; \begin{equation}\label{MK1}\b(x)\cdot(y\cdot z)=(x\cdot y)\cdot \b(z).\end{equation}
  \item The multiplication is Hom-distributive over the addition on both sides;
  \begin{enumerate}
    \item
    \begin{equation}\label{MK2} \a(x)\cdot(y+z)=x\cdot y+x\cdot z, \end{equation}
    \item \begin{equation}\label{MK3}(y+z)\cdot\a(x)=y\cdot x+ z\cdot x,\end{equation}
  \end{enumerate}
  for all $x,~y,~z \in A.$
\end{enumerate}
If $(A,+,\cdot,0,\a,\b)$ admits a unit element $1\in A$ satisfying the following properties
\begin{enumerate}
\item \begin{equation}\label{MK4}x\cdot 1=1\cdot x=\b(x),~\forall x\in A,\end{equation}
\item \begin{equation}\label{MK5}\a(1)=\b(1)=1,\end{equation}
\end{enumerate}
the Hom-ring $A$ is said to be unitary.
\end{df}
\begin{rem}Note that if $\a=id$, we automatically have the definition of Hom-ring in \cite{FG}. \\
When $\a=\b$, a Hom-ring of type $(1)$ should be reffered to as "$\a$-Hom-ring of type $(1)$".
\end{rem}
\begin{pro}
Let $(A,+,\cdot,0,\a,\b)$ be a tuple, where $(A,+,\cdot,0)$ is an abelian Hom-group and $\a$ is a multiplicative map. If $(A,+,\cdot,0,\a,\b)$ satisfies (\ref{MK1}), (\ref{MK2}) or (\ref{MK3}), (\ref{MK4}) and (\ref{MK5}) , then $\b$ is an endomorphism of the abelian  hom-group $(A,+,0,\a)$ and a multiplicative map.
\end{pro}
\begin{prf}

Suppose that $(A,+,\cdot,0,\a,\b)$ satisfies (\ref{MK1}), (\ref{MK2}), (\ref{MK4}) and (\ref{MK5}), so for any $x,~y\in A$, we have
\begin{align*}
  \b(x+y)=1\cdot(x+y)=\a(1)\cdot(x+y)=1\cdot x+1\cdot y=
   \b(x)+\b(y),
\end{align*}
and
\begin{align*}
  \a\circ\b(x)=\a(\b(x))=
   \a(1\cdot x)=
   \a(1)\cdot\a(x)=
   1\cdot\a(x)=
   \b(\a(x)).
\end{align*}
The multiplicativity of $\b$ is obtained  from (\ref{MK1})
\begin{align*}
  \b(x\cdot y)=(x\cdot y)\cdot 1=
   (x\cdot y)\cdot \b(1)=
   \b(x)\cdot(y\cdot 1)=
   \b(x)\b(y).
\end{align*}
\end{prf}
\begin{rem}
As a consequence of the previous proposition, a unitary Hom-ring of type $(1)$ can be defined as a tuple $(A,+,\cdot,0,\a,\b)$ such that $(A,+,0,\a)$ is an abelian Hom-group, $\a$ is a multiplicative map and $(A,+,\cdot,0,\a,\b)$ satisfies the identities (\ref{MK1}), (\ref{MK2}),  (\ref{MK3}), (\ref{MK4}) and (\ref{MK5}).
\end{rem}

Now, we are going to introduce another type of Hom-ring as follows
\begin{df}\label{homring2}
A Hom-ring of type $(2)$ is a  tuple $(A,+,\cdot,0,\a,\b)$ consisting of  a set $A$ together with two binary operations $+:A\times A\rightarrow A$  (the addition) and  $\cdot:A\times A\rightarrow A$ (the multiplication) and two set maps  $\a,\b :A\rightarrow A$, such that:
\begin{enumerate}
  \item $(A,+,0,\a)$ is an abelian Hom-group.
  \item $\b$ is an endomorphism of the abelian  Hom-group $(A,+,0,\a)$, i.e., $\b(x+y)=\b(x)+\b(y),\;\forall x,~y\in A,$ and $\a\circ\b=\b\circ\a.$
  \item $\a$ and $\b$ are multiplicative maps, i.e., $\a(xy)=\a(x)\a(y)$ and  $\b(xy)=\b(x)\b(y)$, for all $x,~y\in A.$
  \item The maps $\a,~\b$ and the product $\cdot$ satisfy $\b(\a^2(x))\cdot(\b(\a(y))\cdot \b^2(z))=(\a^2(x)\cdot \b(\a(y)))\cdot \b^2(\a(z)).$
  \item The multiplication is Hom-distributive over the addition on both sides;
  \begin{enumerate}
    \item $\a^2(x)\cdot\b(y+z)=\a(x)\cdot \b(y)+\a(x)\cdot \b(z)$,
    \item $\a(y+z)\cdot\a(\b(x))=\a(y)\cdot \b(x)+\a(z)\cdot \b(x),$
  \end{enumerate}
  for all $x,~y,~z \in A.$
\end{enumerate}
If $(A,+,\cdot,0,\a,\b)$ admits a unit element $1\in A$ satisfying the following identities
\begin{enumerate}
\item $x\cdot 1=\b(x)$ and $1\cdot x=\a(x),~\forall x\in A,$
\item $\a(1)=\b(1)=1,$
\end{enumerate}
the Hom-ring $A$ is said to be unitary.
\end{df}
\begin{df}\label{homring4}
Let  $(A,+,\cdot,0,\a,\b )$ be a Hom-ring of type $(2)$. It is said to be an "$\a$-Hom-ring of type $(2)$" if $\a=\b.$
\end{df}
\begin{df}
A Hom-ring $(A,+,\cdot,0,\a,\b )$ is called a commutative Hom-ring if for all $x,~y\in A$ we have $x\cdot y=y\cdot x$. Also, a Hom-ring is said to be regular if the twist maps $\a$ and $\b$ are  bijective. The inverse of an element $x$ in a Hom-ring $(A,+,\cdot,0,\a,\b)$ with respect to the addition is named the additive inverse and it is  denoted by $(-x)$.
\end{df}
Now we have introduced these two types of Hom-ring, it's natural to seek for relations among them and provide a way for constructing them.
\begin{lem}\label{lemmering1}
Let $(A,+,\cdot,0,\a,\b)$ be a Hom-ring of type $(1)$.  One has for all $x,~y,~z\in A$:
\begin{enumerate}
  \item $\b(\a^2(x))\cdot(\b(\a(y))\cdot \b^2(z))=\a(\a(x)\cdot\b(y))\cdot\b^3(z).$
  \item $(\a^2(x)\cdot \b(\a(y)))\cdot \b^2(\a(z))=\b\circ\a(\a(x)\cdot(y\cdot z)).$
  \item $\a^2(x)\cdot\b(y+z)=\a(x)\cdot \b(y)+\a(x)\cdot \b(z).$
  \item $\a(y+z)\cdot\a(\b(x))=\a(y)\cdot \b(x)+\a(z)\cdot \b(x).$
\end{enumerate}
\end{lem}
\begin{lem}\label{lemmering2}
Let $(A,+,\cdot,0,\a,\b)$ be a regular Hom-ring of type $(2)$. Then, for all $x,~y,~z\in A$:
\begin{enumerate}
  \item $\b(x)\cdot(y\cdot z)=(x\cdot y)\cdot\a(z).$
  \item $(x\cdot y)\cdot\b(z)=\b(x)\cdot(y\b(\a^{-1}(z)))$.
  \item $\a(x)\cdot(y+z)=x\cdot y+x\cdot z$.
  \item $(y+z)\cdot\a(x)=y\cdot x+ z\cdot x.$
\end{enumerate}
\end{lem}
The desired equivalence is then a simple corollary of this two lemmas:
\begin{pro}
Regular $\a$-Hom-rings of type $(1)$ and $(2)$ are equivalent.
\end{pro}
The following  provides a source of examples for Hom-ring;
\begin{pro}
Let $(A,+,\cdot,0)$ be a ring and $\a,\b:A\rightarrow A$ two commuting homomorphism of  of rings. Defining a new addition $\widetilde{+}:A\times A\rightarrow A$ and a new multiplication $\widetilde{\cdot}:A\times A\rightarrow A$ given by
\begin{itemize}
  \item $x\widetilde{+}y =\a(x+y)=\a(x)+\a(y),$
  \item $x\widetilde{\cdot}y =\b(x\cdot y)=\b(x)\cdot\b(y),$ for all $x,~y\in A$,
\end{itemize}
then $(A,\widetilde{+},\widetilde{\cdot},0,\a,\b)$ is a Hom-ring of type $(1).$ We denote this by $A_{\a,\b}$ and it is called the twist ring of $A.$
\end{pro}
\begin{prf}
The proof is straightforward and left to the reader.
\end{prf}
\begin{pro}
Let $(A,+,\cdot,0,\a,\b)$ be a regular Hom-ring  of type $(1)$. Defining a new addition $+':A\times A\rightarrow A$ and a new multiplication $\cdot':A\times A\rightarrow A$ given by
\begin{itemize}
  \item $x+'y =\a^{-1}(x+y)=\a^{-1}(x)+\a^{-1}(y),$
  \item $x\cdot'y =\b^{-1}(x\cdot y)=\b^{-1}(x)\cdot\b^{-1}(y),$ for all $x,~y\in A$,
\end{itemize}
then $(A,+',\cdot',0)$ is a ring. It is called the compatible ring of the Hom-ring $(A,+,\cdot,0,\a,\b)$.
\end{pro}
\begin{pro}
Let $(A,+,\cdot,0)$ be a ring and $\a,\b:A\rightarrow A$ two commuting homomorphism of rings. Defining a new addition $\widetilde{+}:A\times A\rightarrow A$ and a new multiplication $\widetilde{\cdot}:A\times A\rightarrow A$ given by
\begin{itemize}
  \item $x\widetilde{+}y =\a(x+y)=\a(x)+\a(y),$
  \item $x\widetilde{\cdot}y =\b(x)\cdot \a(y)$, for all $x,~y\in A$,
\end{itemize}
then $(A,\widetilde{+},\widetilde{\cdot},0,\a,\b)$ is a Hom-ring of type $(2).$ We denote this by $A_{\a,\b}$ and it is called the twist ring of $A.$
\end{pro}
\begin{pro}
Let $(A,+,\cdot,0,\a,\b)$ be a regular Hom-ring  of type $(2)$. Defining a new addition $+':A\times A\rightarrow A$ and a new multiplication $\cdot':A\times A\rightarrow A$ given by
\begin{itemize}
  \item $x+'y =\a^{-1}(x+y)=\a^{-1}(x)+\a^{-1}(y),$
  \item $x\cdot'y =\b^{-1}(x)\cdot\a^{-1}(y),$ for all $x,~y\in A$,
\end{itemize}
then $(A,+',\cdot',0)$ is a ring, which is called the compatible ring of the Hom-ring $(A,+,\cdot,0,\a,\b)$.
\end{pro}
In what follows we  provide some  properties and define  Hom-ring homomorphisms and Hom-subrings.
\begin{pro}\label{R2}
Let $(A,+,\cdot,0,\a,\b )$ be a unitary regular hom-ring of type $(1)$ or of type $(2).$
\begin{enumerate}
\item For all $a\in A$ , we have $a\cdot 0=0\cdot a=0.$
\item Let $a,~b\in A$. Then , the additive inverse of $(a\cdot b)$ is $(-a)\cdot b=a\cdot(-b)=-(a\cdot b)$. One can notes  that $\a(-a)=-\a(a).$
\end{enumerate}
\end{pro}
\begin{prf}\label{invrs}
We have $$0\cdot \a(a)=(0+0)\cdot\a(a)=0\cdot a+0\cdot a.$$
Let $c$ be the additive inverse of $0\cdot a$, i.e.,  $c+0\cdot a=0$. Then we have
\begin{align*}
0=\a(c+0\cdot a)=
   \a(c)+0\cdot\a(a)=
   \a(c)+(0\cdot a+0\cdot a))=
  (c+0\cdot a)+\a(0\cdot a)=
  \a^2(0\cdot a).
\end{align*}
Since $\a$ is injective, then $0\cdot a=0$. The second property is a consequence of the first property.
\end{prf}

\begin{df}
Let $(A,+_A,\cdot_A,0_A,\a_A,\b_A)$ and $(B,+_B,\cdot_B,0_B,\a_B,\b_B)$ be two unitary hom-rings.
\begin{itemize}
  \item  A homomorphism of Hom-rings is a map $f: A\rightarrow B$ such that
   \begin{enumerate}
     \item $f(0_A)=0_B$ and $f(1_A)=1_B$
     \item $f(x+_A y)=f(x)+_B f(y)$
     \item $f(x\cdot_A y)=f(x)\cdot_B f(y).$
     \item $f\circ\a_A=\a_B\circ f$, and $f\circ\b_A=\b_B\circ f.$
   \end{enumerate}
   \item The map $f$ is called a weak homomorphism if it satisfies the second and the third condition.
  \item The Hom-rings $(A,+_A,\cdot_A,0_A,\a_A,\b_A)$ and $(B,+_B,\cdot_B,0_B,\a_B,\b_B)$ are called isomorphic if $f$ is bijective.
  \item If $f$ is the inclusion map, we say that  $(A,+,\cdot,0_A,\a,\b)$ is a Hom-subring of $(B,+,\cdot,0_B,\a,\b).$
\end{itemize}
\end{df}
\begin{rem}
Let $(A,+,\cdot,0_A,\a,\b)$ be  a Hom-subring of $(B,+,\cdot,0_B,\a,\b)$. Then for all $x\in A$, we have (obviously) $\a(x)\in A$ and $\b(x)\in A$  . Also, a Hom-subring of a Hom-ring $(A,+,\cdot,0_A,\a,\b)$ is a Hom-subgroup of $(A,+,0_A,\a)$ stable by the multiplication and the map $\b$.
 \end{rem}
\begin{df}
Let $(A,+,\cdot,0,\a,\b)$ be a unitary regular Hom-ring of type $(1)$ (resp. regular Hom-ring of type $(2)$). The center $Z(A)$ (resp. center $Z'(A)$) of $A$ is the set of all $x\in A$ where $x\cdot a=a\cdot x$ (resp. $x\cdot\b(a)=\a(a)\cdot x$), for all $a\in A$.
\end{df}
\begin{pro}
Let $(A,+,\cdot,0,\a)$ be a unitary regular Hom-ring.
\begin{enumerate}
  \item If $(A,+,\cdot,0,\a)$ is of type $(1)$, then $Z(A)$ is a Hom-subring of $A$.
  \item If $(A,+,\cdot,0,\a)$ is of type $(2)$, then $Z'(A)$ is a Hom-subring of $A$.
\end{enumerate}
\end{pro}
\begin{prf}
\begin{enumerate}
\item Let $(A,+,\cdot,0,\a)$ be a Hom-ring of type $(1)$
We have, for all $a\in A$, $0\cdot a=a\cdot 0=0$ and $1\cdot a=a\cdot 1=\b(a)$. Then $0$ and $1$ are in $Z(A).$
Let $x,~y\in Z(A)$. Then for all $a\in A$, we have
\begin{enumerate}
  \item $(x+y)\cdot a=(x+y)\cdot\a(\a^{-1}(a))=x\cdot\a^{-1}(a)+y\cdot\a^{-1}(a)=\a^{-1}(a)\cdot x+\a^{-1}(a)\cdot y=a\cdot(x+y).$
  Thus $Z(A)$ is  closed under addition  of $A$.
  \item If $x\in Z(A)$ and $a\in A$ then $-(x\cdot a)=-(a\cdot x)$. Furthermore, we have $-(x\cdot a)=(-x)\cdot a$ and
  $-(a\cdot x) =a\cdot(-x).$  So $-x\in Z(A)$.
  \item $(x\cdot y)\cdot a=(x\cdot y)\cdot\a(\a^{-1}(a))=\a(x)\cdot(y\cdot\a^{-1}(a))=\a(x)\cdot(\a^{-1}(a)\cdot y)=(x\cdot\a^{-1}(a))\cdot\a(y)=(\a^{-1}(a)\cdot x)\cdot\a(y)=a\cdot(x\cdot y).$
      Thus $Z(A)$ is closed under multiplication of $A$.
\end{enumerate}
\item Let $(A,+,\cdot,0,\a)$ be a Hom-ring of type $(2)$.
We have, for all $a\in A$, $0\cdot \b(a)=\a(a)\cdot 0=0$ and $1\cdot \b(a)=\a(a)\cdot 1=\b\circ\a(a)$. Then $0$ and $1$ are in $Z'(A).$
Let $x,~y\in Z'(A)$. Then for all $a\in A$, we have
\begin{enumerate}
  \item $(x+y)\cdot \b(a)=\a^{-1}(\a^{-1}(x)+\a^{-1}(y))\cdot\b(\a(\a^{-1}(a)))=x\cdot\b(\a^{-1}(a))+y\cdot\b(\a^{-1}(a))=a\cdot x+a\cdot y=\a(a)\cdot(x+y).$
  Thus $Z'(A)$ is  closed under addition  of $A$.
  \item If $x\in Z'(A)$ and $a\in A$ then $-(x\cdot a)=-(a\cdot x)$. Furthermore, we have $-(x\cdot a)=(-x)\cdot a$ and
  $-(a\cdot x) =a\cdot(-x).$ Then $-x\in Z'(A).$
  \item $(x\cdot y)\cdot \b(a)=\b(x)\cdot(y\cdot\b(\a^{-1}(a)))=\b(x)\cdot(a\cdot y)=(x\cdot a)\cdot\a(y)=\a(a)\cdot(x\cdot y).$
      Thus $Z'(A)$ is closed under multiplication of $A$.
\end{enumerate}
\end{enumerate}
\end{prf}
\subsection{Examples}
\begin{pro}
Let $(A,+,\a_A)$ be a regular abelian Hom-group and  consider the group homomorphisms from $A$ into $A$. Then $(End(A),+,\circ,0,\a_{End(A)},id)$  is a Hom-ring of type $(1)$ such that
\begin{enumerate}
 \item The addition of two such homomorphisms may be defined pointwise to produce another group homomorphism. Explicitly, given two such homomorphisms $f$ and $g$, the sum of $f$ and $g$ is the homomorphism $(f+g)(x)=f(x)+g(x)$,
 \item The composition $\circ$ of two such homomorphisms $f$ and $g$ is explicitly $(f\circ g)(x)=f(g(x))$,
 \item $\a_{End(A)} (f)=\a_A\circ f$.
 \end{enumerate}
\end{pro}
\begin{prf}
The proof is straightforward and left to the reader.
\end{prf}
{\bf Polynomial Hom-Rings:} Let $A$ be a ring and $R=A[X_1,\cdot\cdot\cdot,X_n]$ be the polynomial ring in $n$ variables.\\ Let us recall that any element of $R$ is uniquely written as a linear combination, with coefficients in $A$, of monomes $X^n=X_1^{m_1} ... X_n^{m_n}$, where $m$ runs over the set of sequences of $n$ nonnegative integers.\\
   Take  a ring endomorphism  of $A$, $\a:A\rightarrow A\subset R.$ Then a ring endomorphism $\widetilde{\a}$ of R is uniquely determined by the
$n$ polynomials $\widetilde{\a}(X_i)=\sum_r a_{i,r}X^r=\sum_r a_{i,r}X_1^{r_1} ... X_n^{r_n}$ for $1\leq i\leq n.$ \\Let
$\widetilde{\a}(P)= \sum_m\a(a_{m})\widetilde{\a}(X_1)^{m_1} ... \widetilde{\a}(X_n)^{m_n}, $ for $P=\sum_m a_{m}X_1^{m_1} ... X_n^{m_n}\in R.$\\
For  $P=\sum_m a_{m}X_1^{m_1} ... X_n^{m_n}$ and $Q=\sum_m b_{m}X_1^{m_1} ... X_n^{m_n}$, define $\widehat{+}$ by
$$P\widehat{+} Q=\sum_m \a(a_{m}+b_m)\widetilde{\a}(X_1)^{m_1} ... \widetilde{\a}(X_n)^{m_n}, $$
and for  $P=\sum_m a_{m}X_1^{m_1} ... X_n^{m_n}$ and $Q=\sum_l b_{l}X_1^{l_1} ... X_n^{l_n}$, define $\widehat{\cdot}$ by
$$P\widehat{\cdot} Q=\sum_d\sum_{m+l=d} \a(a_{m})\a(b_l)\widetilde{\a}(X_1)^{d_1} ... \widetilde{\a}(X_n)^{d_n}.$$
Then $R_{\widetilde{\a}}=(R,\widehat{+},\widehat{\cdot},\widetilde{\a})$ is a $\a-$Hom-ring of type $(1)$.

{\bf Twisted Group-Ring:} Let $G$ be a group and $\K G$ the corresponding group-ring over $\K$. the set $\K G$ is generated by $\{e_g : g\in G\}$.If $\a:G\rightarrow G$ is a group homomorphism, then it can be extended to a ring endomorphism of $\K G$ by setting
      $$\a\Big(\sum_{g\in G}a_g e_g\Big)=\sum_{g\in G}a_g \a(e_g)=\sum_{g\in G}a_g e_{\a(g)}.$$
      Consider the usual ring structure on $\K G$. Then, we define a  $\a-$Hom-ring of type $(1),~(\K G,\widehat{+},\widehat{\cdot},\widetilde{\a}) $ over $\K G$ by setting:
      \begin{align*}
       \sum_{g\in G}a_g e_g\widehat{+}\sum_{g\in G}b_g e_g=\sum_{g\in G}(a_g+b_g) e_{\a(g)},\ \ \ (\sum_{g\in G}a_g e_g)\widehat{\cdot }(\sum_{h\in G}b_h e_h)=\sum_{gh\in G}a_gb_g e_{\a(gh)}.
      \end{align*}

\section{Modules over $\a$-Hom-rings of type $(1)$}
We introduce in the following the structure of module over $\a$-Hom-ring of type $(1).$
The properties of modules for other types Hom-ring will be discussed in a forthcoming paper.\\
In the rest of this article, the word Hom-ring means an $\a$-Hom-ring of type $(1).$
\subsection{Definitions and Properties}
\begin{df}
Let $(A,+_A,\cdot,0,\a)$ be a unitary  Hom-ring with unit $1_A$. A left $A$-module $M$ consists of an abelian Hom-group $(M,+_M,\beta)$ and an operation $\lambda_l: A\otimes M\rightarrow M$  such that for all $a,~b $ in $A$ and $m_1,~m_2$ in $M$, we have:
\begin{equation}\label{M1}
\lambda_l(\a(a)\otimes(m_1+_M m_2))=\lambda_l(a\otimes m_1)+_M\lambda_l(a\otimes m_1),
\end{equation}
\begin{equation}\label{M2}
 \lambda_l((a+_A b)\otimes \b_M(m_1))=\lambda_l(a\otimes m_1)+_M\lambda_l(b\otimes m_1),
\end{equation}
\begin{equation}\label{M3}
 \lambda_l((a\cdot b)\otimes \b(m_1))=\lambda_l(\a(a)\otimes \lambda_l(b\otimes m),
\end{equation}
\begin{equation}\label{M4}
  \lambda_l(1_A\otimes m_1)=\b(m_1).
\end{equation}
The operation of the Hom-ring on $M$ is called scalar multiplication, and it is usually written as $am$ for $a$ in $A$ and $m$ in $M$. The notation $_AM$ indicates a left $A$-module $M$. If $\b$ is invertible,  such left module $_AM$ is said regular. A right $A$-module $M$ or $M_A$ is defined similarly, except that the ring acts on the right; i.e., the operation takes the form $\ld_r:M\otimes A\rightarrow M$, and the above axioms are written with a scalar multiplication  on the right.
\end{df}
\begin{ex}
  If $A$ is a  unitary  Hom-ring, then we may consider $A$ as either a left or right A-module, where the action is given by the multiplication.
\end{ex}
\begin{pro}
Let $(A,+_A,\cdot,0,\a)$ be a unitary  Hom-ring. If $A$ is commutative, then  these left modules are the same as right modules.
\end{pro}
\begin{prf}
  Let $M$ be a right $A$-module, i.e., there exists $\ld_r:M\otimes A\longrightarrow M$ by $\ld_r(m\otimes a)=ma$. Now, we define $\ld_l:A\otimes M\longrightarrow M$ by $\ld_l(a\otimes m)=ma$. Using the commutativity of the Hom-ring $A$, we can verify that $\ld_l$ defines a  left module structure for $M$ over $A.$
\end{prf}
\begin{df}
 If $R$ and $S$ are two Hom-rings, then an $R-S$-bimodule is an abelian Hom- group $M$  such that:
 \begin{enumerate}
   \item $M$ is a left $R$-module and a right $S$-module,
   \item $(rm)\a_S(s)=\a_R(r)(ms)$.
 \end{enumerate}
\end{df}
\begin{df}
A bimodule is both a left $A$-module with action $\ld_l$ and
a right $A$-module with action $\ld_r$ satisfying the compatibility condition
\begin{equation}\label{eq1}
  \ld_r\circ(\ld_l\otimes\a)=\ld_l\circ(\a\otimes\ld_r).
\end{equation}
It is  simply called $A$-module.
\end{df}
\begin{pro}
Let $(A,+_A,\cdot,0_A,\a)$ be a unitary  Hom-ring and $(M,+_M,0_M,\b)$ be a regular left $A$-module.\begin{enumerate}
\item  For all $m$ in $M$ , we have  $0_Am=0_M.$
\item If moreover $a\in A$ admits an additive inverse of index $r\leq n$, then the additive inverse of $am$ on $M$ is $-(am)=(-a)m=a(-m)$.  We can observe that $\b(-m)=-\b(m)$.
    \end{enumerate}
\end{pro}
\begin{prf}
We can show that the first property in the same way as the first property of Proposition \ref{invrs}. For the second property, we have
$$\a^{n-1}(a)\b^n(m)+\a^{n-1}(-a)\b^n(m)=\a^n(a+(-a))\b^{n+1}(m)=0\b^{n-1}(m)=0.$$
In the same way we check the other identity.
\end{prf}
\begin{rem}
If $M$ is a  module over a Hom-ring $A$, then two  equations (\ref{M3}) and (\ref{M4}) give  that \begin{equation}\label{M5}
\b(am)=\a(a)\b(m).\end{equation}
One can get, for all $m$ element in right (regular) $A$-module $M$, $m0_A=0_M$.
\end{rem}

\subsection{Module Constructions }
\begin{pro}\label{CSTR1}
  Let $(A,+_A,\cdot_A,0,\a)$ be a unitary regular Hom-ring, $(M,+_M,\b)$ a regular $A$-module (left), $(A,+_A',\cdot')$ the compatible ring. If $\begin{array}{l|rcl}
  \ld_l:&A&\longrightarrow& End(M)\\
        &a&\longmapsto& \ld_l(a),
        \end{array}$ is the action of the Hom-ring $A$ on $M$, then $(M,+_M'=\b^{-1}\circ +_M)$ is a (left) $(A,+_A',\cdot')$-module with respect the action $\ld_l'=\b^{-1}\ld_l.$
\end{pro}

\begin{prf}
  For the action $\ld_l'$, the conditions for $(M,+_M')$ to be an $(A,+_A',\cdot')$-module is written in the form
  \begin{equation}\label{MM1}
    \b^{-1}(a(m_1+_Mm_2))=\b^{-2}(am_1+_Mam_2),
  \end{equation}
\begin{equation}\label{MM2}
    \b^{-1}((a+_Ab)m_1))=\b^{-2}(am_1+_Mbm_1),
  \end{equation}
 \begin{equation}\label{MM3}
    \b^{-1}(a\b^{-1}(bm_1)=\b^{-1}(\a^{-1}(a\cdot_Ab)m_1),
  \end{equation}
  \begin{equation}\label{MM4}
    \b^{-1}(1_Am_1)=m_1.
  \end{equation}
Also, one can check that  equation (\ref{M5}) is equivalent to \begin{equation}\label{MM5}
\ld_l\a(a)=\b\ld_l(a)\b^{-1}.\end{equation}
Then the result can  easily follow.
\end{prf}
From Proposition \ref{CSTR1}, we can get a method of building (left) regular modules  of a unitary regular Hom-ring:
\begin{thm}\label{modmod}
  Let $(A,+_A,\cdot_A,\a)$ be a unitary regular Hom-ring with the compatible ring $(A,+_A',\cdot_A')$.
\begin{enumerate}
  \item If $(M,+_M,\b)$ is a (left) regular $A$-module then \begin{equation}\label{cstr2}
\b(a\rhd m)=\a(a)\rhd\b(m),~\;\forall a\in A,~\;m\in M,
  \end{equation}
  where $a\rhd m=\b^{-1}(am)$ defines a $(A,+_A',\cdot_A')$-module structure for $(M,\b^{-1}\circ +_M)$, it is called the compatible $A$-module of $(M,+_M,\b)$.
  \item Suppose that $(M,+_M')$ is a (left) module over $(A,+_A',\cdot_A')$ for the action $a\rhd m$. If it exists an automorphism of the group $(M,+_M')$ such that (\ref{cstr2}) is satisfied, then $(M,+_M=\b\circ +_M',\b)$ is a (left) module over $(A,+_A,\cdot_A,\a)$ with respect to the action given by  $am=\b(a\rhd m)$, for all $a\in A$ and $m\in M$.
\end{enumerate}
\end{thm}
\begin{prf}\item
  \begin{enumerate}
    \item $(M,+_M,\b)$ is a (left) regular $A$-module, i.e., $\b$ is invertible. Then from (\ref{M5}), we can get easily (\ref{cstr2}).
    \item Let $a,~b\in A $ and $m_1,~m_2\in M$, from (\ref{cstr2}) and the properties of $(M,+_M')$ as a(left) $(A,+_A',\cdot_A')$-module, we have
    \begin{align}
     & \a(a)(m_1+_Mm_2) =\b(\a(a)\rhd\b(m_1+_M'm_2))\nonumber \\
       &=
       \a^2(a)\rhd(\b^2(m_1)+_M'\b^2(m_2))
       =
       \a^2(a)\rhd\b^2(m_1)+_M'\a^2(a)\rhd \b^2(m_2)\nonumber \\
       &=
       \b(\a(a)\rhd\b(m_1))+_M'\b(\a(a)\rhd \b(m_2))
       =
       \a(a)\rhd\b(m_1))+_M\b(\a(a)\rhd \b(m_2)\nonumber\\
       &=
       \b(a\rhd m_1)+_M\b(a\rhd m_2)\nonumber\\
       &=am_1+_Mam_2,
    \end{align}
    \begin{align}
      &(a+_Ab)\b(m_1) =\b((a+_Ab)\rhd\b(m_1))
       =
       \a^2((a+_A'b))\rhd\b^2(m_1)\nonumber \\
       &=
       \a^2(a)\rhd\b^2(m_1)+_M'\a^2(b)\rhd \b^2(m_1)
       =
       \b(\a(a)\rhd\b(m_1))+_M'\b(\a(b)\rhd \b(m_1))\nonumber  \\
       &=
       \a(a)\rhd\b(m_1))+_M\b(\a(b)\rhd \b(m_1)
       =
       \b(a\rhd m_1)+_M\b(b\rhd m_1)\nonumber\\
       &=am_1+_Mbm_1,
    \end{align}
      \begin{align}
     & \a(a)(bm_1) =\b(\a(a)\rhd(\a(b)\rhd\b(m)))
       =
       \a^2(a)\rhd\b(\a(b)\rhd\b(m_1))=
       \a^2(a)\rhd(\a^2(b)\rhd\b^2(m_1))\nonumber \\
       &=
       (\a^2(a)\cdot'\a^2(b))\rhd\b^2(m_1)=
       \a(a\cdot b)\rhd\b(m_1))=
       \b((a\cdot b)\rhd m_1)\nonumber\\
       &=
       (a\cdot b)\b(m),
    \end{align}
    \begin{align}
    1_Am=\b(1_A\rhd m_1)=
    \a(1_A)\rhd\b(m_1)=
    1_A\rhd \b(m_1)=\b(m_1).
    \end{align}
  \end{enumerate}
Therefore $(M,+_M=\b\circ +_M',\b)$ is a (left) module over $(A,+_A,\cdot_A,\a)$ with respect the action defined by  $am=\b(a\rhd m)$, for all $a\in A$ and $m\in M.$
\end{prf}

\subsection{Tensor product of modules}
\begin{df}
Let $(R,\a_R)$ be a Hom-ring. Let $(A,\alpha_A) $ be a right $R$-module, $(B,\alpha_B)$ a left $R$-module and $(C,\alpha_C)$ be an abelian Hom-group (all abelian Hom-groups having the law of composition denoted by $+$). A map $f:A\times B\rightarrow C$ from the Cartesian product of  $A$ and $B$ into  $C$ is called Hom-R-Bilinear  if
\begin{enumerate}
  \item $f(a_1+a_2,\a_B(b))=f(a_1,b)+f(a_2,b)$,
  \item $f(\a_A(a),b_1+b_2)=f(a,b_1)+f(a,b_2)$,
  \item $f(\a_A(a),\a_B(b))=\a_C\circ f(a,b)$,
  \item $f(ar,\a_B(b))=f(\a_A(a),rb)$,
\end{enumerate}
for all $a,~a_1,~a_2\in A,~b,~b_1,~b_2\in B$ and $r\in R$. The set of all such Hom-R-bilinear maps from $A\times B$ to $C$ is denoted by $HL_R(A,B;C).$
\end{df}
\begin{df}
Let $(R,\a_R)$ be a Hom-ring, $(A,\alpha_A) $ be a right regular $R$-module and $(B,\alpha_B)$ be a left regular $R$-module. The tensor product of $A$ and $B$, which is denoted by $A\otimes_R B$, is a regular abelian Hom-group together with a Hom-R-Bilinear map
$$\t:A\times B\rightarrow A\otimes_R B,$$
such that the following universal property holds:
for any Hom-R-bilinear map $f:A\times_R B \rightarrow C$, there exists a unique Hom-module homomorphism $\widetilde{f}:A\otimes_R B \rightarrow C$ such that the diagram:\\
 \vskip 0.3cm
$$\xymatrix{ A\times B \ar[rr]^\t \ar^{f}[rd] && A\otimes_R B \ar[ld]^{\widetilde{f}} \\ & C },$$
commutes, that is $\widetilde{f}\circ \t=f.$
\end{df}
\begin{pro}
Let $R,~S$ be Hom-rings. If $A$ is a regular $R-S$-bimodule and $B$ is a left regular $S$-module, then $A\otimes_S B$ is a left $R$-module satisfying $$r(a\otimes b)=(ra)\otimes b.$$
\end{pro}
\begin{prf}
  For each $s\in S$, we have
  \begin{align*}
    \a_R(r)(as)\otimes \a_B(b) &= (ra)\a_S(s)\otimes \a_B(b) 
     =
      \a_A(ra)\otimes (\a_S(s)n)
     =
     (\a_R(r)\a_A(a)) \otimes (\a_S(s)n).
  \end{align*}
\end{prf}

\section{Semisimple Modules and  Hom-rings}
\begin{df}
Let $A$ be a Hom-ring and $(M,\b_1)$ and $(N,\b_2)$ be two (left) $A$-modules. A  homomrophism of Hom-groups $f:N\rightarrow M$ is said to be a homomorphism of $A$-modules if $f(n)=\a_A(a)f(n)$, for all $a\in A$ and $n$ in $M.$ \\
If $f$ is the inclusion map, we say that $N$ is an $A$-submodule of $M.$
\end{df}
\begin{rem}
An $A$-submodule of $M$ is a Hom-subgroup $N$ of $M $ which is itself a module under the action of Hom-ring elements.
\end{rem}
Note that, for any $A$-module $M$, both $\{0\}$ and $M$ are trivial $A$-submodules of $M$. Let us now recall that the direct sum of two Hom-groups $(H,\b_1)$ and $(K,\b_2)$ (\cite{H2}) is the set
$$H\oplus K=\{(h,k)| h\in H,~k\in K\},$$
together with the component-wise operations of the groups $H$ and $K$ (that is, if $(h_1,k_1),~(h_2,k_2)\in H\oplus K$, then $(h_1,k_1)(h_2,k_2)=(h_1h_2,k_1k_2) $) and the map $\b_1\oplus\b_2$ such as $\b_1\oplus\b_2 (h,k)=(\b_1(h),\b_2(k)),~\forall (h,k)\in H\oplus K$.
There is a similar structure with modules.
\begin{df}
  Let $(M,\b_1)$ and $(N,\b_2)$ be modules over a Hom-ring $A$.  We define the direct sum of
$M$ and $N$ by
\begin{equation}\label{mod1}
M\bigoplus N=\{m\oplus n| m\in M,~n\in M\}.
\end{equation}
The structure of $A$-module of $M\bigoplus N$ is given by
\begin{equation}\label{mod2}
  (m_1\oplus n_1)\oplus(m_2\oplus n_2)=(m_1+_Mm_2)\oplus(n_1+_Nn_2),
\end{equation}
\begin{equation}\label{mod3}
  (\b_1\oplus\b_2)(m\oplus n)=\b_1(m)\oplus\b_2(n),
\end{equation}
\begin{equation}\label{mod4}
  a(m\oplus n)=am\oplus an,\;\;\mbox{where}\;\;a\in A.
\end{equation}
\end{df}

\begin{df}
Let $A$ be a unitary (regular)  Hom-ring. An $A$-module $M$ is said to be simple if it is non-zero, and has no submodules except $\{0\}$ and itself. If $M$ is a direct sum of certain simple submodules, $M$ is called semisimple.
\end{df}

\begin{lem}\label{simple1}\item
 Let $(M,+_M,\b)$ be a (left) module over a unitary Hom-ring $A$. Then $(ker\b,+_M,\b)$ is an $A$-submodule.
\end{lem}
\begin{prf}
  For all $m\in Ker\b$, by $\b(am)=\a(a)\b(m)=0,$ for all $a\in A$, it is easy to show the lemma.
\end{prf}
\begin{thm}\label{simple2}\item
Let $(M,+_M,\b)$ be a (left)simple module over a unitary Hom-ring $A$. Then $M$ is regular, i.e., the map $\b$ is invertible.
\end{thm}
\begin{prf}
According to Lemma \ref{simple1}, $(ker\b,+_M,\b)$ is an $A$-submodule. Since the module $M$ is simple, either $ker\b=\{0\}$ or $ker\b=M$. The $A$-module $M$ is simple, then  $ker\b=\{0\}$. So $\b$ is invertible.
\end{prf}
\begin{pro}
Let $(A,+_A,\cdot_A,\a)$ be a unitary regular Hom-ring with the compatible ring $(A,+_A',\cdot_A')$. Suppose that $(M,+_M')$ is a simple (left) module over $(A,+_A',\cdot_A')$ for the action $a\rhd m$. If it exists an automorphism of the group $(M,+_M')$ such that (\ref{cstr2}) is satisfied, then $(M,+_M=\b\circ +_M',\b)$ is a simple (left) module over $(A,+_A,\cdot_A,\a)$ with respect to the action given by  $am=\b(a\rhd m)$, for all $a\in A$ and $m\in M$.
\end{pro}
\begin{prf}
  Indeed, let $N$ be a non-zero $A$-submodule of $(M,+_M=\b\circ +_M',\b)$, that is $N\subset M$ and $(N,+_M=\b\circ +_M',\b)$ is an $A$-module with respect to the action defined by $a\otimes n\mapsto an,~\forall a\in A,~n\in N$. According to Theorem \ref{modmod}, $(N,+_M')$ is a module over the ring $(A,+_A',\cdot_A')$ with respect to the action $a\otimes n\mapsto a\rhd n=\b^{-1}(an)$, where $\rhd$ defines the action of the ring $(A,+_A',\cdot_A')$ on $(M,+_M')$. Since $(M,+_M')$ is simple, then $N=M$ as groups with respect to the addition $+_M'$. So  from the property $\b(N)=N$, we can show the result.
\end{prf}
The following theorem gives a characterization of (left) modules of a regular unitary Hom-ring.
\begin{thm}\label{modmod1}
Let $(A,+_A,\cdot_A,\a)$ be a unitary regular Hom-ring with the compatible ring $(A,+_A',\cdot_A')$ and $(M, +_M,\b)$ be a simple regular  $A$-module. Then the compatible  group $(M,+_M'=\b^{-1}\circ +_M)$ is a  semisimple module over the ring $(A,+_A',\cdot_A')$.
\end{thm}
\begin{prf}
Suppose that $(M, +_M,\b)$ is a regular  $A$-module  for the action $a\otimes n\mapsto am,~\forall a\in A,~m\in M$, then according to Theorem \ref{modmod}, $(M,+_M'=\b^{-1}\circ +_M)$ is a module over $(A,+_A',\cdot_A')$ for the action $a\otimes n\mapsto a\rhd m=\b^{-1}(am)$. Let $m$ be a non-zero element in $M$ and let $k$ be the smallest integer such that $\b^{k-1}(m)\neq m$ and $\b^{k}(m)= m$. The group of $(M,+_M'=\b^{-1}\circ +_M)$, $N=A\rhd m\bigoplus A\rhd\b(m)\bigoplus\cdot\cdot\cdot\bigoplus A\rhd \b^{k-1}(m)$ is a submodule of $(M,+_M')$ over the ring $(A,+_A',\cdot_A')$. Furthermore,  we have $\b(N)=N=Am\bigoplus A\b(m)\bigoplus\cdot\cdot\cdot\bigoplus A\b^{k-1}(m)$. Consequently $N$ is a submodule of $(M, +_M,\b)$ over $(A,+_A,\cdot_A,\a)$. So the simplicity of $(M, +_M,\b)$ gives $M=N$ as a module over $(M,+_M'=\b^{-1}\circ +_M)$ . Therefore, as a module over $(A,+_A',\cdot_A')$, we have $M=N=A\rhd m\bigoplus A\rhd\b(m)\bigoplus\cdot\cdot\cdot\bigoplus A\rhd \b^{k-1}(m),$ which shows the semisimplicity of the module  $(M,+_M'=\b^{-1}\circ +_M)$ over $(A,+_A',\cdot_A')$ for the action $a\otimes n\mapsto a\rhd m=\b^{-1}(am)$.
\end{prf}
\begin{df}
A Hom-ring $A$ is called semisimple (resp. simple) if $A$  is semisimple (resp. simple) as an $A$-module.
\end{df}
\begin{thm}
  Simples unitary Hom-rings are regular Hom-rings.
\end{thm}
\begin{prf}
The proof is straightforward.
\end{prf}
\begin{thm}\label{modmod1}
Let $(A,+_A,\cdot_A,\a)$ be a unitary  Hom-ring with the compatible ring $(A,+_A',\cdot_A')$. Then the compatible  ring $(A,+_A',\cdot_A')$ is   semisimple.
\end{thm}
\begin{prf}
The proof is straightforward.
\end{prf}
\bibliographystyle{amsplain}

\end{document}